\documentclass{amsart}
\usepackage{amsthm,amstext,amsmath,amscd,amssymb,latexsym}
\usepackage{color}
\usepackage[matrix,arrow]{xy}

\usepackage[left=3cm,top=2.75cm,right=3cm]{geometry}

\input xy
\xyoption{all}
\usepackage{amsfonts,enumerate,array,hyperref}

\newcommand{\A}{{\mathbb A}}

\newcommand{\C}{{\mathbb C}}

\newcommand{\K}{{\mathbb K}}

\newcommand{\N}{{\mathbb N}}
\renewcommand{\P}{{\mathbb P}}

\newcommand{\T}{{\mathbb T}}

\newcommand{\Z}{{\mathbb Z}}

\newcommand{\ko}{{\mathcal O}}
\newcommand{\cO}{{\mathcal O}}

\newcommand{\s}{\mathcal}

\newcommand{\sI}{{\s I}}

\newcommand{\tZ}{\widetilde{Z}}
\newcommand{\punkt}{\HHspace{-.3ex}\raise.15ex\HHbox to1ex{\HHuge.}}

\newcommand{\paper}{: \begin{it}}
\newcommand{\jour }{, \end{it}}

\newtheorem{theorem}{Theorem}[section]
\newtheorem{lemma}[theorem]{Lemma}
\newtheorem{proposition}[theorem]{Proposition}
\newtheorem{corollary}[theorem]{Corollary}
\newtheorem{conjecture}[theorem]{Conjecture}

\theoremstyle{definition}
\newtheorem{definition}[theorem]{Definition}

\newtheorem{example}[theorem]{Example}

\theoremstyle{remark}
\newtheorem{remark}[theorem]{Remark}

\title[Secant varieties of Segre-Veronese varieties]{On the dimensions
  of secant varieties of Segre-Veronese varieties}
\author{Hirotachi Abo}
\address{Department of Mathematics, University of Idaho, Moscow, ID 83844, USA} \email{abo@uidaho.edu}
\author{Maria Chiara Brambilla}
\address{Dipartimento di Scienze Matematiche,  Universit\`a Politecnica delle Marche, Ancona, Italy}
\email{brambilla@dipmat.univpm.it}


\thanks{
The first author is partly supported by NSF grant DMS-0901816.
The second author is partially supported by Italian MIUR
and is member of GNSAGA-INDAM.}

\begin{document}
\begin{abstract}
This paper explores the dimensions of higher secant varieties to Segre-Veronese
varieties. The main goal of this paper is to introduce two different
inductive techniques.  These techniques enable one to reduce the
computation of the dimension of the secant variety in a high
dimensional case  to the computation of  the dimensions of
secant varieties in low dimensional cases. As an application of these
inductive approaches, we will prove non-defectivity of secant
varieties of certain two-factor Segre-Veronese varieties. We also use
these methods to  give a complete classification of defective
$s^{\mathrm{th}}$  Segre-Veronese varieties for small $s$. In the
final section, we propose a conjecture about defective two-factor
Segre-Veronese varieties.
\end{abstract}
\maketitle
\section{Introduction}
\label{sec:intro}
In many applications, it is natural to represent a collection of data
as a multi-indexed list. Alternatively, one can think of the data as a
multi-dimensional array. A mathematical framework that includes the
study of multi-dimensional arrays is through parameter spaces of
tensors.

Every tensor can be written as a linear combination of so-called {\it
  decomposable tensors}. A tensor is said to have {\it rank} $s$ if it
can be written as a linear combination of $s$ decomposable tensors
(but not fewer).   Note that there are higher rank tensors that
can be written as the limit of lower rank tensors. A tensor is
said to have {\it border rank} $s$ if it can be expressed as the
limit of rank $s$ tensors, but not as the limit of rank $s-1$
tensors. For more details on tensor rank and tensor border rank,
we refer the reader, for example, to \cite{L}. An interesting
question is ``Given a positive integer $s$, what is the dimension
of the parameter space of tensors  with border rank at most $s$?".
In the following few paragraphs, we will formulate this problem as
a classical problem in algebraic geometry.

Let $k$ be a positive integer. For each $i \in \{1, \dots, k\}$, let
$V_i$ be a vector space of dimension $n_i+1$ over $\C$, $n_1\leq
\cdots \leq n_k$. The collection of decomposable tensors can be
``embedded'' into the $N$-dimensional vector space $\bigotimes_{i=1}^k
V_i$, where $N = \prod_{i=1}^k(n_i+1)$.  Projectivizing to account for
the effect of scalars, we have a  {\it Segre map}
$\prod_{i=1}^k\P(V_i) \rightarrow \P\left(\bigotimes_{i=1}^k V_i
\right)$.  The image of this map, denoted $X$,  is called the {\it
  Segre variety}.

A {\it secant $(s-1)$-plane} to $X$ is a linear subspace that passes
through $s$ linearly independent points of $X$. Each point on the
secant $(s-1)$-plane is a linear combination of $s$ points on $X$ and
can be identified with a tensor which is a linear combination of $s$
fixed decomposable tensors. The Zariski closure of the set of all
points which lie on a secant $(s-1)$-plane, i.e., the set of all
tensors that can be written as the sum of $s$ decomposable tensors,
is called the $s^{\mathrm {th}}$ {\it secant variety} of $X$ and
denoted by $\sigma_s(X)$.
The variety $\sigma_s(X)$
parameterizes tensors with border rank at most $s$.  Thus the
aforementioned question is equivalent to the question about ``What is
the dimension of $\sigma_s(X)$?".

Since $\sigma_s(X)\subset\P^{N-1}$ is the closure of the union of secant
$(s-1)$-planes to $X$,  the following inequality holds:
\[
\dim \sigma_s(X) \leq \min \left\{N-1,  s\left(1+\sum_{i=1}^kn_i\right)-1\right\}.
\]
We say that $\sigma_s(X)$ has the {\it expected dimension} if the
equality holds.
The Segre variety $X$ has a {\it defective $s^{\mathrm  {th}}$ secant
  variety} if $\sigma_s(X)$ does not have the expected dimension. In
particular, $X$ is called {\it defective} if $X$ has a defective
$s^{\mathrm{th}}$ secant variety for some $s$. For example, if $k=2$,
then  $X$ corresponds to the parameter space of rank one $(n_1+1)
\times (n_2+1)$ matrices, and the points of $\sigma_s(X)$ correspond
to $(n_1+1) \times (n_2+1)$ matrices that can be written as the sum of
$s$ (or fewer) rank one matrices of the same size.
Thus the affine cone over $\sigma_s(X)$ can be identified with the general determinantal
variety $M_k$ of $n_1 \times n_2$ matrices of rank $s$ or less. Recall that $M_k$ has codimension
$(n_1+1-s)(n_2+1-s)$
(see for example \cite{Har} for more details on determinantal varieties). So if $2 \leq s \leq \min \{n_1, n_2\}$, then the dimension of $\sigma_s(X)$ is strictly smaller
than the expected one. Therefore, most of secant varieties of Segre varieties with two factors are defective.
On the other hand, there are only a
few families of defective Segre varieties known to exist for $k \geq
3$.  It is therefore desirable to classify defective Segre varieties.

There are other categories of tensors such as symmetric tensors,
alternating tensors and mixed regular and symmetric
tensors. Those tensors also arise very naturally throughout physics,
computer science, engineering as well as mathematics.

The concepts of rank and border rank of regular tensors can  be
extended to tensors in other categories. The geometry of
decomposable tensors in each of these categories can be analogously exploited:
Veronese varieties, Grassmann varieties and Segre-Veronese varieties can be thought of as
parameter spaces of decomposable symmetric tensors, decomposable
alternating tensors and decomposable mixed regular and symmetric
tensors respectively, and questions about rank of tensors in each
category are related to questions about secant varieties of the
corresponding varieties.

A well known classification of the defective Veronese varieties was
completed in a series of papers by Alexander and
Hirschowitz~\cite{AH}.  There are corresponding conjecturally complete
lists of defective Segre varieties~\cite{AOP2} and Grassmann
varieties~\cite{BDG}. Defective secant varieties of Segre-Veronese
varieties are, however,  less well-understood, although considerable
efforts have been already made to complete the list of such varieties
(see for example, \cite{CGG2}, \cite{CCh}, \cite{Bal}, \cite{CC},
\cite{O}, \cite{Abr}).  Even the classification of defective
two-factor Segre-Veronese varieties
is still far from complete.

One of the main goals of this paper is to provide several tools to
study secant varieties of Segre-Veronese varieties.
In order to classify defective Segre-Veronese varieties, a crucial
step is to prove the existence of a large family of non-defective such
varieties. A powerful tool to establish non-defectivity of large
classes of Segre-Veronese varieties is the inductive approach based on
specialization techniques, which consist in placing a certain number
of points on a chosen divisor. For a given $\mathbf{n}=(n_1, \dots,
n_k)  \in \N^k$, we denote $\P^{n_1} \times \cdots \times \P^{n_k}$ by $\P^{\mathbf{n}}$. Let $X_{\mathbf{n}}^{\mathbf{a}}$ be the
Segre-Veronese variety obtained by embedding $\P^{\mathbf{n}}$ in
$\P^{\prod_{i=1}^k {n_i+a_i \choose a_i}-1}$  by the morphism given by
$\cO(\mathbf{a})$ with $\mathbf{a}=(a_1, \dots, a_k) \in \N^k$.
Thanks to the classical theorem called Terracini's lemma (see Theorem \ref{th:terracini} for a more detailed statement of Terracini's lemma)  ,
it is easy to see that the problem of determining the dimension of
$\sigma_s(X_{\mathbf{n}}^{\mathbf{a}})$ is equivalent to the problem
of determining  the value of the Hilbert function
$h_{\P^{\mathbf{n}}}(Z, \cdot)$ of a collection $Z$ of $s$ general
double points in $\P^{\mathbf{n}}$ at $\mathbf{a}$, i.e.,
\[
h_{\P^{\mathbf{n}}}(Z, \mathbf{a}) = \dim H^0 ( \P^{\mathbf{n}},
\cO(\mathbf{a})) - \dim H^0 ( \P^{\mathbf{n}}, \sI_Z(\mathbf{a})).
\]
Suppose that $a_{1} \geq 2$. Denote by $\mathbf{n'}$ and $\mathbf{a}'$
the $k$-tuples $(n_1-1, n_2, \dots, n_k)$ and $(a_1-1, a_2, \dots,
a_k)$ respectively. Given a $\P^{\mathbf{n}'} \subset
\P^{\mathbf{n}}$,
we have a short exact sequence
\[
 0 \rightarrow \sI_{\widetilde{Z}}(\mathbf{a}') \rightarrow \sI_Z(\mathbf{a}) \rightarrow
 \sI_{Z \cap \P^{\mathbf{n'}},  \P^{\mathbf{n'}}}(\mathbf{a}) \rightarrow 0,
\]
where $\widetilde{Z}$ is the residual scheme of $Z$ with respect to
$\P^{\mathbf{n}'}$ and $Z \cap \P^{\mathbf{n}'}$ is the trace of $Z$ on the hyperplane. This short exact sequence gives rise to the
so-called {\it Castelnuovo inequality}
\[
h_{\P^{\mathbf{n}}}(Z, \mathbf{a}) \geq h_{\P^{\mathbf{n}}}(\widetilde{Z}, \mathbf{a}')
+ h_{\P^{\mathbf{n}'}}(Z \cap \P^{\mathbf{n}'}, \mathbf{a}).
\]
Thus, we can conclude that
\begin{itemize}
\item[(a)] if $h_{\P^{\mathbf{n}}}(\widetilde{Z}, \mathbf{a}')$ and
$h_{\P^{\mathbf{n}'}}(Z \cap \P^{\mathbf{n}'}, \mathbf{a}')$ are the
expected values and
\item[(b)] if the degrees of $\widetilde{Z}$ and $Z\cap
  \P^{\mathbf{n}'}$ are both less than or both greater than $\dim
  H^0(\P^{\mathbf{n}}, \cO(\mathbf{a}'))$ and $\dim
  H^0(\P^{\mathbf{n}'}, \cO(\mathbf{a}))$ respectively,
\end{itemize}
then $h_{\P^{\mathbf{n}}}(Z, \mathbf{a})$ is also the expected value.
By semicontinuity,  the Hilbert function of
a general collection of $s$ double points in $\P^{\mathbf{n}}$
has the expected value at $\mathbf{a}$.

The problem is, however, that it may or may not be possible to arrange that Condition (b) is satisfied.
In Section~\ref{sec:induction} we  generalize the {\it m\'ethode
  d'Horace diff\'erentielle} of Alexander and Hirschowitz~\cite{AH} to give a way around this
numerical obstacle.
The precise statement of our version of the  Horace method can be found in
Theorem~\ref{th:horace}. For the reader's convenience,  we state the same theorem in a slightly different format than Theorem~\ref{th:horace} below.
\begin{theorem}
 \label{th:horaceIntro}
Let $a_1 \geq 3$.
Let $\mathbf{n}'=(n_1-1, n_2, \dots, n_k)$, let $\mathbf{a}'=(a_1-1, a_2, \dots, a_k)$, and let
$\mathbf{a}''=(a_1-2, a_2, \dots, a_k)$.
For a given positive integer $s$, let $s'$ and $\epsilon$ be the quotient and remainder
when dividing $s\left(1+\sum_{i=1}^k n_i\right)-{n_1+a_1-1 \choose a_1-1}\prod_{i=2}^k{n_i+a_i \choose a_i}$ by $\sum_{i=1}^k
n_i$. Suppose that $s'\geq \varepsilon$.
If $\sigma_{s'}(X_{\mathbf{n'}, \mathbf{a}})$, $\sigma_{s-s'}(X_{\mathbf{n}, \mathbf{a'}})$, and
$\sigma_{s-s'-\epsilon}(X_{\mathbf{n}, \mathbf{a''}})$ have the expected dimension and if
\begin{eqnarray}
\label{eq:superabundant}
(s-s'-\epsilon)\left(1+\sum_{i=1}^k n_i\right) \geq {n_1+a_1-2 \choose a_1-2}\prod_{i=2}^k {n_i+a_i \choose a_i},
\end{eqnarray}
then $\sigma_{s}(X_{\mathbf{n}, \mathbf{a}})$ also has the expected dimension.
\end{theorem}
This theorem enables one to check whether or not $\sigma_s(X_{\mathbf{n}, \mathbf{a}})$
has the expected dimension by induction on $\mathbf{n}$ and $\mathbf{a}$.
It cannot however be applied to $\sigma_s(X_{\mathbf{n}, \mathbf{a}})$
if $\mathbf{a}$ is small. The theorem requires that one of the $a_i$'s
is at least $3$, so one cannot use it when every $a_i$ is less than or
equal to two. In addition, if at least one of the degrees is $1$, it
is frequent that Inequality~(\ref{eq:superabundant}) does not hold. In
Section~\ref{sec:induction}, we therefore develop a different
inductive approach for computing the dimensions of secant varieties of
such Segre-Veronese varieties. This approach allows one to place a
certain number of points not only on a hypersurface, but also on a
subvariety (see Theorem~\ref{th:segrein} for a more precise
statement). Note that a similar approach was successfully applied to study
secant varieties of Segre varieties in~\cite{AOP2}.

In order to apply these inductive approaches, we need some initial
cases regarding either dimensions or degrees. The class of secant
varieties of two-factor Segre-Veronese varieties
can be viewed as one of such initial cases. In
Section~\ref{sec:two-factor},  we will study secant varieties of such
Segre-Veronese varieties. The main goal of this section is to prove
the following theorem:
\begin{theorem}
\label{th:n1abintro}
Let $n,a\ge1$, $b\ge3$, $\mathbf{n} = (n,1)$ and $\mathbf{a} = (a,b)$.
Then $X_{\mathbf{n},\mathbf{a}}$ is not defective except if $(n,a,b)=(n,2,2k)$.
\end{theorem}

\noindent
We will restate and prove this theorem in Section~\ref{sec:two-factor} (see Corollary~\ref{coro3.14}).

The strength of Theorem~\ref{th:horaceIntro}  is to reduce establishing the existence of a large number of families
of non-defective Segre-Veronese varieties to establishing the existence of only  a small number of
families of non-defective cases. We will prove the following theorem as an application of Theorems~\ref{th:horaceIntro} and \ref{th:n1abintro} to demonstrate the power of Theorem~\ref{th:horaceIntro}:
\begin{theorem}
\label{th:boundintro}
Suppose that $X_{\mathbf{n},\mathbf{a}}$ is not defective for every
$\mathbf{n}$ and for $\mathbf{a}=(3,3)$, $(3,4)$ and $(4,4)$. Then
$X_{\mathbf{n},\mathbf{a}}$ is not defective for every $\mathbf{n}$
and for every $\mathbf{a}=(a,b)$ such that $a, b\geq 3$.
\end{theorem}

\noindent
This theorem will also be restated and proved in Section~\ref{sec:two-factor} (see Theorem~\ref{th:bound}).

As we shall see in Section~\ref{sec:induction}, using a randomized
algorithm which employs  Terracini's lemma, we can compute the
dimension of  $\sigma_s(X_{\mathbf{n},\mathbf{a}})$  for a given $s
\in \N$ and for given $\mathbf{n}, \mathbf{a} \in \N^k$. Based on our
experiments using this randomized algorithm, we expect
 that there
are no defective Segre-Veronese varieties $X_{\mathbf{n},\mathbf{a}}$
for any
$\mathbf{n}$ if $\mathbf{a}=(3,3)$, $(3,4)$ or $(4,4)$. Thus
Theorem~\ref{th:boundintro}  suggests the following conjecture:
\begin{conjecture}
 Let $\mathbf{n}$ and $\mathbf{a}$ be pairs of positive integers.  If
 $\mathbf{a} \geq (3,3)$, there are no defective two-factor
 Segre-Veronese varieties $X_{\mathbf{n}, \mathbf{a}}$ for all
 $\mathbf{n} \in \mathbb{N}^2$.
 \end{conjecture}
In Section~\ref{sec:classification}, we apply the inductive procedures
developed in Section~\ref{sec:induction}  to classify all
the defective $s^{\mathrm{th}}$ secant varieties of Segre-Veronese
varieties for each $s \in \{2,3,4\}$.

Section~\ref{sec:conjecture} provides a conjecturally
complete list of defective secant varieties of two-factor
Segre-Veronese varieties.  In addition to evidence provided by our theorems,
further evidence in support of the
conjecture was obtained via the computational experiments we carried
out with {\tt Macaulay2},  a computer algebra system developed by Dan Grayson and Mike
Stillman~\cite{GS}.

\section{Inductive techniques}
\label{sec:induction}
For each $i \in \{1, \dots, k\}$, let $V_i$ be a $(n_i+1)$-dimensional
vector space over $\C$ and let $\P^{n_i}=\P(V_i)$.
Given two $k$-tuples $\mathbf{n}=(n_1,\cdots, n_k)$ and
$\mathbf{m}=(m_1,\cdots, m_k)$,
we write $\mathbf{n}\le \mathbf{m}$ when $n_i\le m_i$ for all $i$.
Unless otherwise stated,
$\mathbf{n}$, $\mathbf{n}'$, $\mathbf{a}$, $\mathbf{a}'$ and  $\mathbf{a}''$
denote
$(n_1,\cdots, n_k)$, $(n_1-1,n_2, \cdots, n_k)$, $(a_1, \cdots, a_k)$,
$(a_1-1, a_2,\cdots, a_k)$ and  $(a_1-2, a_2,\cdots, a_k) \in \N^k$ respectively.
We write  $\P^{\mathbf{n}}$ for $\prod_{i=1}^k \P^{n_i}$ and
$X_{\mathbf{n},\mathbf{a}}$ for the Segre-Veronese variety embedded in
$\P^{N-1}$ by $\ko_{\P^{\mathbf{n}}}(\mathbf{a})$, where
$N=\prod_{i=1}^k {n_i+a_i \choose a_i}$.
Let $N_R={n_1+a_1-1 \choose a_1-1}\prod_{i=2}^k {n_i+a_i \choose a_i}$
and $N_T={n_1+a_1-1 \choose a_1}\prod_{i=2}^k {n_i+a_i \choose
  a_i}$.
Let $R=\C\left[x_{0,1}, \dots, x_{n_1,1}, \dots ,x_{0,k}, \dots,
  x_{n_k,k}\right]$ and note that it can be thought of as an
$\N^k$-graded ring in the obvious way.

Let $\sigma_s(X_{\mathbf{n}, \mathbf{a}})$ be the $s^{\mathrm{th}}$
secant variety of $X_{\mathbf{n}, \mathbf{a}}$, i.e., the Zariski
closure of the union of linear subspaces spanned by $s$-tuples of
points on $X_{\mathbf{n}, \mathbf{a}}$.
We now explain how to translate the problem of computing the
dimension of $\sigma_s(X_{\mathbf{n}, \mathbf{a}})$
into a question about 
the value of the Hilbert function of the ideal of $s$ double points on
$\P^{\mathbf{n}}$ at $\mathbf{a}$.
Let
$\T_p(X_{\mathbf{n}, \mathbf{a}})$ be the projective tangent
space to $X_{\mathbf{n}, \mathbf{a}}$ at a point $p$.
The following well known result describes the tangent space of
$\sigma_s(X_{\mathbf{n}, \mathbf{a}})$:
\begin{theorem}[Terracini's lemma]
Let $p_1, \dots, p_s$  be generic points of $X_{\mathbf{n},
  \mathbf{a}}$ and let $q$ be a generic point of $\langle p_1, \dots,
p_s \rangle$. Then
\[
 \T_q [\sigma_s(X_{\mathbf{n}, \mathbf{a}})] =
\left\langle \T_{p_1}(X_{\mathbf{n}, \mathbf{a}}),\ldots, \T_{p_s}(X_{\mathbf{n}, \mathbf{a}})\right\rangle,
\]
where $\T_q[\sigma_s(X_{\mathbf{n}, \mathbf{a}})]$ is the projective
tangent space to $\sigma_s(X_{\mathbf{n}, \mathbf{a}})$ at $q \in
\sigma_s(X_{\mathbf{n}, \mathbf{a}})$.
\label{th:terracini}
\end{theorem}
\begin{remark}
\label{rem:computation}
Let $\mathbf{n}$ and $\mathbf{a}$ be $k$-tuples of non-negative integers.
Let $k$ be a positive integer. For an $i \in \{1, \dots, k\}$, let
$V_i$ be an $(n_i+1)$-dimensional vector space over $\C$ and let $v_i
\in V_i \setminus \{0\}$. Denote by $p \in X_{\mathbf{n},\mathbf{a}}$
the equivalence class containing $v_1^{a_1} \otimes \cdots \otimes
v_k^{a_k}$.
Then the affine cone  over $\T_p(X_{\mathbf{n},  \mathbf{a}})$ in
$\bigotimes_{i=1}^k S_{a_i} V_i$ is
\[
C[\T_p(X_{\mathbf{n},  \mathbf{a}})]
=\sum_{i=1}^k v_1^{a_1} \otimes \cdots \otimes v_i^{a_i-1}V_i \otimes \cdots \otimes v_k^{a_k}.
\]  In particular, $C[\T_p(X_{\mathbf{n},  \mathbf{a}})]$ can be
represented by a $\left[\sum_{i=1}^k (n_i+1)\right] \times N$ matrix
$A_p$.
Thus Terracini's lemma can be used to estimate the dimension of
$\sigma_s(X_{\mathbf{n},  \mathbf{a}})$ as follows: First choose
randomly $s$ points $p_1, \dots, p_s$ on $X_{\mathbf{n},
  \mathbf{a}}$. Next, compute the matrix representation $A_{p_i}$ for
each $C[\T_{p_i}(X_{\mathbf{n},  \mathbf{a}})]$.  Let $A$ be the
matrix $\left(\begin{array}{c} A_1\\\vdots \\A_n\end{array}\right)$.
It follows from
Terracini's lemma that $\dim \sigma_s(X_{\mathbf{n},  \mathbf{a}})
\geq \mathrm{rank}(A)-1$. By semi-continuity, the equality holds if $
\mathrm{rank}(A) = \min\left\{s\left(1+\sum_{i=1}^k n_i\right), N\right\}$, because
$\dim  \sigma_s(X_{\mathbf{n},  \mathbf{a}})  \leq
\min\left\{s\left(1+\sum_{i=1}^k n_i\right)-1, N-1\right\}$.  Finally, we would like to
stress that although \linebreak
$\mathrm{rank}(A)\not= \min\left\{1+\sum_{i=1}^k
  n_i, N\right\}$ is a strong evidence that $\sigma_s(X_{\mathbf{n},
  \mathbf{a}})$ is defective, it cannot be used to prove defectivity.
\end{remark}

Note that $H^0(\P^{\mathbf{n}},
\mathcal{O}_{\P^{\mathbf{n}}}(\mathbf{a}))$ can be identified with the
set of hyperplanes in $\P^N$. Since the condition that a hyperplane $H
\subset \P^N$ contains $\T_p(X_{\mathbf{n}, \mathbf{a}})$ is
equivalent to the condition that $H\cap X_{\mathbf{n}, \mathbf{a}}$ contains
the first infinitesimal neighborhood of $p$, the
elements of $H^0(\P^{\mathbf{n}}, \mathcal{I}_p^2(\mathbf{a}))$ can be
viewed as hyperplanes containing $\T_p(X_{\mathbf{n}, \mathbf{a}})$.
Let $Z$ be a collection of $s$ double points on $\P^{\mathbf{n}}$
and let $\sI_Z$ be its ideal sheaf.
Terracini's lemma implies that $\dim \sigma_s(X_{\mathbf{n}, \mathbf{a}})$ is equal to the value of the Hilbert function $h_{\P^{\mathbf{n}}}(Z,\cdot)$ of $Z$ at $\mathbf{a}$.
Hence proving that  $\sigma_s(X_{\mathbf{n}, \mathbf{a}})$ has the expected
dimension is equivalent to proving that
\[
h_{\P^{\mathbf{n}}}(Z,\mathbf{a})= \min \left\{
 s\left(1+\sum_{i=1}^k n_i\right), \ N \right\}.
\]
The following definition is analogous to Definition 3.2 in \cite{AOP2}:
\begin{definition}
\label{abundance}
Let $\mathbf{n}, \mathbf{a} \in \N^k$, let $s$ be a non-negative
integer and $Z$ a zero-dimensional subscheme of $\P^n$.
A triple $(\mathbf{n}; \mathbf{a};Z)$ is said to be {\it subabundant}
(resp. {\it superabundant}) if $\deg Z \leq N$ (resp. $\deg Z \geq N$).
The triple $(\mathbf{n}; \mathbf{a};Z)$ is said to be {\it equiabundant} if it is both subabundant and superabundant.
We say that two triples {\it have the same abundancy} if both of them are either superabundant, or subabundant.
We say that $T(\mathbf{n}; \mathbf{a};Z)$ is {\it true} if
$h_{\P^{\mathbf{n}}}(Z,.)$ has the expected value at $\mathbf{a}$.
If $Z$ is a collection of $s$ general double points, we write $T(\mathbf{n};
\mathbf{a};s)$ instead of $T(\mathbf{n}; \mathbf{a};Z)$ and
$(\mathbf{n}; \mathbf{a};s)$ instead of $(\mathbf{n}; \mathbf{a};Z)$.
We say that $T(\mathbf{n}; \mathbf{a})$ is true if
$T(\mathbf{n}; \mathbf{a};s)$ is true for every $s \geq 0$.
\end{definition}
Assume that $a_1\ge2$. Let $H$ be a hypersurface defined by a  linear form in
$R_{(1,0,\ldots,0)}$.
For a given zero-dimensional subscheme $Z$, we denote by $\tZ$ the
{\it residual} of $Z$ with respect to $H$,
i.e.\ the subscheme whose ideal is $\sI_Z:\sI_H$.
The scheme $Z\cap H$ is called the {\it trace} of $Z$.
>From the restriction exact sequence
\[
0\rightarrow \sI_{\widetilde{Z}}(\mathbf{a'})\rightarrow
\sI_Z(\mathbf{a})\rightarrow \sI_{Z\cap H}(\mathbf{a})\rightarrow 0,
\]
we easily get  the so-called Castelnuovo inequality
$$h_{\P^{\mathbf{n}}}(Z,\mathbf{a})\geq
h_{\P^{\mathbf{n}}}(\tZ,\mathbf{a'})+
h_{\P^{\mathbf{n'}}}(Z\cap H,\mathbf{a}).$$
>From this inequality it is easy to prove the following basic Horace lemma:
\begin{theorem}
\label{th:horace-easy}
Let $a_1\ge2$, let $Z$ be a zero-dimensional subscheme of $\P^{\mathbf{n}}$ and
let $H$ be a hyperplane defined by a linear form in $R_{(1,0,\ldots,0)}$.
\begin{itemize}
\item[(i)] If  $h_{\P^{\mathbf{n}}}(\tZ,\mathbf{a'})$ and
$h_{\P^{\mathbf{n'}}}(Z\cap H,\mathbf{a})$ are equal to the expected value;
\item[(ii)] if $(\mathbf{n}'; \mathbf{a}; Z\cap H)$ and
$(\mathbf{n};\mathbf{a}' ; \tZ)$ have the same abundancy,
\end{itemize}
then $h_{\P^{\mathbf{n}}}(Z,\mathbf{a})$ is also the expected value.
\end{theorem}
\begin{lemma}\label{chandler}
Assume that $a_1\ge2$.
Let $Z$ be a subscheme of $\P^{\mathbf{n}}$ and
let $H$ a hyperplane defined by a linear form in $R_{(1,0,\ldots,0)}$.
Then there exists a collection $\Phi$ of $u$ general points in $H$ such that
\[
h_{\P^{\mathbf{n}}}(Z\cup\Phi, \mathbf{a})
= h_{\P^{\mathbf{n}}}(Z, \mathbf{a})+u
\]
if and only if $u$ satisfies
\begin{equation}\label{condition-H}
h_{\P^{\mathbf{n}}}(Z, \mathbf{a})+u\le
h_{\P^{\mathbf{n}}}(\tZ, \mathbf{a'})+\binom{n_1-1+a_1}{n_1-1}
\prod_{i=2}^k {n_i+a_i \choose a_i}.
\end{equation}
\end{lemma}
\begin{proof}
This lemma is an easy generalization of Lemma~$3$
in~\cite{Ch}. One can prove our statement exactly in the same way as
in~\cite{Ch}, and thus we omit the proof.
\end{proof}
In the following example, we show how to combine
Theorem~\ref{th:horace-easy} with Lemma~\ref{chandler}, in order to reduce
computing the dimension of the secant variety of a Segre-Veronese
variety to computing the dimensions of secant varieties of smaller
Segre-Veronese varieties.
 \begin{example}
\label{example-horace-easy}
Let $\mathbf{n} = (1,1)$ and let $\mathbf{a} = (3,3)$.  Let $p_1,
\dots, p_5 \in \P^{\mathbf{n}}$ and let $Z = \{p_1^2, \dots,
p_5^2\}$. Specialize two points,
say $p_4$ and $p_5$, to $H = \P^0 \times \P^1 \subset (\P^1)^2$.  Then
$\widetilde{Z}$ consists of three double points and two simple points;
while $Z \cap H$ consists of two double points in $H$. So  both
$(1,1;2,3;\widetilde{Z})$ and $(0,1;3,3;Z \cap H) = (1;3;Z \cap H) = (1;3;2)$
are subabundant.  It is well known that $T(1;3;2)$ is true.
We therefore want to prove the truth of  $T(1,1;2,3;\widetilde{Z})$.

Note that  the inequality
\begin{eqnarray*}
11 &= & 9+2 \\
&=& h_{\P^{\mathbf{n}}}(\{p_1^2,p_2^2,p_3^2\}, (2,3))+2 \\
 &\le & h_{\P^{\mathbf{n}}}(\{p_1^2,p_2^2,p_3^2\}, (1,3))+
{1+3 \choose 3} \\
&= &8+4=12,
\end{eqnarray*}
holds. Thus, by Lemma~\ref{chandler},
the expected value of the Hilbert function of $\widetilde{Z}$ at $(2,3)$ is
\[
h_{\P^{\mathbf{n}}}(\widetilde{Z},(2,3)) =
h_{\P^{\mathbf{n}}}(\{p_1^2,p_2^2,p_3^2\},(2,3))+2 = 11.
\]
Additionally, Theorem 2.1 in \cite{CGG2}
implies that  $T(1,1;2,3;3)$ and $T(1,1;1,3;3)$ are true.
Thus, \linebreak
$T(1,1;2,3;\widetilde{Z})$ is true. Therefore,
the truth of $T(\mathbf{n};\mathbf{a};5)$ follows from Theorem~\ref{th:horace-easy}.
\end{example}
As already stated in Section~\ref{sec:intro}, one cannot always arrange  that Condition (ii) in
Theorem~\ref{th:horace-easy} is satisfied. We illustrate it in the following example:
\begin{example}
\label{counter}
Let  $\mathbf{n}=(2,2)$, let $\mathbf{a}=(4,4)$, let $p_1, \dots,
p_{45} \in \P^{\mathbf{n}}$ and let $Z=\{p_1^2, \dots, p_{45}^2\}$. To
prove the truth of $T(\mathbf{n};\mathbf{a};s)$, we want to specialize
a certain number of points among the $p_i$'s, say $p_1, \dots, p_{s'}$,
to $H \simeq \P^1 \times \P^2 \subset (\P^2)^2$ in such a way that
$(1,2; \mathbf{a}; s')$ and $(\mathbf{n};3,4 ; \tZ)$ have the same
abundancy. This means that  they must be equiabundant, because
$(\mathbf{n};\mathbf{a};45)$ is equiabundant.  It is not possible,
however,  to find such an integer $s'$, because ${1+4 \choose 4}{2+4
  \choose 4}/(1+2+1) \not\in \Z$. Thus one cannot apply
Theorem~\ref{th:horace-easy} to show that $T(\mathbf{n};\mathbf{a};s)$
is true.
\end{example}
One of the main goals of this section is to generalize the differential
Horace method introduced by Alexander and Hirschowitz to Segre-Veronese varieties in order to side step numerical obstacles like above.

Given a linear system $\mathcal{D}$ on $\P^{\mathbf{n}}$,
we say that a scheme $Z$ is {\it $\mathcal{D}$-independent} if the value
$h_{\P^{\mathbf{n}}}(Z, \mathcal{D})= \dim H^0(\P^{\mathbf{n}}, \mathcal{D})-
\dim H^0(\P^{\mathbf{n}}, \sI_Z\otimes\mathcal{D})$ equals the degree
of the scheme $Z$.
The following lemma is also due to Chandler (see  \cite[Lemma 6.1]{BO} for a detailed proof):
\begin{lemma}
\label{curvilinear}
Let $Z\subset \P^{\mathbf{n}}$ be a zero-dimensional scheme
contained in a finite collection of double points and let $\mathcal{D}$ be a
linear system on $\P^{\mathbf{n}}$.
Then $Z$ is $\mathcal{D}$-independent if and only if
every curvilinear subscheme $\zeta$ of $Z$ is $\mathcal{D}$-independent.
\end{lemma}
We are now able to prove the {\it m\'ethode d'Horace diff\'erentielle} for Segre-Veronese varieties.
\begin{theorem}
\label{th:horace}
Let $a_1 \geq 3$.
For a given non-negative integer  $s$, let  $s'$ and $\varepsilon$ be
the quotient and remainder in the division of $s\left(1+\sum_{i=1}^k
  n_i\right)-N_R$ by  $\sum_{i=1}^k n_i$. Suppose that $s'\ge\varepsilon$.
If $T(\mathbf{n}'; \mathbf{a}; s')$, $T(\mathbf{n};\mathbf{a}' ; s-s')$
and $T(\mathbf{n};\mathbf{a}''; s-s'-\varepsilon)$ are all true and if
$(\mathbf{n};\mathbf{a}''; s-s'-\varepsilon)$ is superabundant, then
$T(\mathbf{n}; \mathbf{a};s) $ is also true.
\end{theorem}
\begin{proof}
Here we only focus on the case when  $(\mathbf{n};
\mathbf{a};s)$ is subabundant, because the remaining case can be
proved in a similar manner.

\vspace{2mm}
\noindent
{\sc Step 1}.
By assumption,
$N_R=\left(1+\sum_{i=1}^k n_i\right)(s-s')-\varepsilon+s'$, and since $s'\ge\varepsilon$
we have that $(\mathbf{n}; \mathbf{a}'; s-s')$ is subabundant.
This implies that since $T(\mathbf{n}; \mathbf{a}'; s-s')$ holds by assumption, then the Hilbert function
$h_{\P^{\mathbf{n}}}(Z,\mathbf{a}')$ has the expected value for any subscheme $Z$ of a collection of $s-s'$ general
double points.

Now choose a hyperplane $H$ defined  by  a linear form in $R_{(1, 0, \dots, 0)}$.
Let $\Gamma=\{\gamma^1,\ldots ,\gamma^{\varepsilon}\}$ be a collection of $\varepsilon$
general points contained in $H$ and $\Sigma$ a collection of $s-s'-\varepsilon$ points
not contained in $H$. Let $Z=\Gamma^2_{|H}\cup\Sigma^2$.
Then from what we say above it follows
\[
h_{\P^{\mathbf{n}}}(Z,\mathbf{a}')=\min\left\{\left(1+\sum_{i=1}^k
    n_i\right)(s-s')-\varepsilon, N_R\right\}=\left(1+\sum_{i=1}^k n_i\right)(s-s')-\varepsilon.
\]

\vspace{2mm}
\noindent
{\sc Step 2.}
Now we want to add to $Z$ a collection $\Phi$ of $s'$ simple points contained in $H$ in such a way that
\begin{eqnarray}
\label{eq:chandlereq}
h_{\P^{\mathbf{n}}}(Z\cup\Phi,\mathbf{a}')=h_{\P^{\mathbf{n}}}(Z,\mathbf{a}')+s'.
\end{eqnarray}
By Lemma~\ref{chandler} we can do this if
\[h_{\P^{\mathbf{n}}}(Z,\mathbf{a}')+s'\le
h_{\P^{\mathbf{n}}}(\Sigma^2, \mathbf{a}'')+{n_1+a_1-2 \choose
  a_1-1}\prod_{i=2}^k {n_i+a_i \choose a_i}.\]
By assumption, $T(\mathbf{n}; \mathbf{a}'';s-s'-\varepsilon)$ is true
and $(\mathbf{n}; \mathbf{a}'';s-s'-\varepsilon)$ is superabundant,
which implies
\[
h_{\P^{\mathbf{n}}}(\Sigma^2, \mathbf{a}'')+{n_1+a_1-2 \choose
  a_1-1}\prod_{i=2}^k {n_i+a_i \choose a_i} =N_R.
\]
On the other hand, by Step 1 we know that
$h_{\P^{\mathbf{n}}}(Z,\mathbf{a}')+s'=N_{R}$, then
Equality~(\ref{eq:chandlereq}) follows.

\vspace{2mm}
\noindent
{\sc Step 3.}
>From the assumption that $(\mathbf{n}, \mathbf{a},s)$ is subabundant  and the definition of $s'$ and $\varepsilon$ it follows that
$$s'\left(\sum_{i=1}^k n_i\right)+\varepsilon =s\left(\sum_{i=1}^k n_i+1\right)- N_R\le N-N_R=N_T.$$
Since $T(\mathbf{n}', \mathbf{a},s')$ holds by assumption,
the scheme $(\Gamma\cup\Phi^2_{|H})\subset H$ has Hilbert function
$$h_{\P^{\mathbf{n'}}}(\Gamma
\cup\Phi^2_{|H},\mathbf{a})=s'\sum_{i=1}^k n_i+\varepsilon$$

\vspace{2mm}

Now, for $(t_{1},\ldots,t_{\varepsilon})\in\K^{\varepsilon}$,  choose a flat family
of general points
$\Delta_{(t_1,\ldots, t_\varepsilon)}=\{\delta_{t_{1}}^{1},\ldots,\delta_{t_{\varepsilon}}^{\varepsilon}\}\subseteq \P^{\mathbf{n}}$ and a
family of hyperplanes $\{H_{t_{1}},\ldots,H_{t_{\varepsilon}}\}$ defined by linear forms in $R_{(1,0,\ldots,0)}$ such that

$\bullet\quad \delta_{t_{i}}^{i}\in H_{t_{i}}$  for any $t_{i}$, and
any $i=1,\ldots,\varepsilon$,

$\bullet\quad\delta_{t_{i}}^{i}\not\in H$ for any $t_{i}\neq0$,  and
any $i=1,\ldots,\varepsilon$,

$\bullet\quad H_0= H$ and $\delta_0^{i}=\gamma^{i}\in H$,  for any $i=1,\ldots,\varepsilon$.

\noindent Now let us consider the following schemes:

$\bullet\quad\Delta_{(t_1,\ldots, t_\varepsilon)}^2=\{\delta^1_{t_1},\dots,\delta_{t_\varepsilon}^{\varepsilon}\}^2$,
  notice that
$\Delta_{(0,\ldots, 0)}^2={\Gamma}^2$;

$\bullet\quad {\Phi}^2$, where $\Phi$ is the collection of the $s'$ points introduced in Step $2$;

$\bullet\quad {\Sigma}^2$, the collection of the $s-s'-\varepsilon$ double points introduced in Step $1$.

In order to prove $T(\mathbf{n}; \mathbf{a};s)$ it is enough to prove
the following claim.

\vspace{2mm}
\noindent
{\sc Claim.}
There exists $(t_1,\ldots, t_\varepsilon)$ such that the scheme
$\Delta_{(t_1,\ldots, t_\varepsilon)}^2$ is independent with respect
to the linear system $\sI_{\Phi^2\cup\Sigma^2}\otimes
\ko_{\P^{\mathbf{n}}}(\mathbf{a})$.

\vspace{2mm}
\noindent
{\it Proof of the claim.}
Assume that the claim is false. Then by Lemma
\ref{curvilinear} for all $(t_1,\ldots, t_\varepsilon)$ there exist pairs
$(\delta_{t_i}^i,\eta_{t_i}^i)$ for $i=1,\ldots,\varepsilon$, with $\eta_{t_i}^i$
a curvilinear scheme supported in $\delta_{t_i}^i$ (hence the length of $\eta_{t_i}^i$ is $2$ for all $i$) and contained in
$\Delta_{(t_1,\ldots, t_\varepsilon)}^2$ such that
\begin{equation}\label{assurdo}
h_{{\P}^n}(\Phi^2\cup\Sigma^2\cup\eta_{t_1}^1\cup\ldots \cup \eta_{t_\varepsilon}^{\varepsilon},\mathbf{a})<
\left(1+\sum_{i=1}^k n_i\right)(s-\varepsilon )+2\varepsilon.
\end{equation}
Let $\eta_0^i$ be the limit of $\eta_{t_i}^i$, for
$i=1,\ldots,\varepsilon$.
Suppose that $\eta_0^i\not\subset H$ for $i\in
F\subseteq\{1,\ldots,\varepsilon\}$ and $\eta_0^i\subset H$
for $i\in G=\{1,\ldots,\varepsilon\}\setminus F$.
Given $t\in\K$, let us denote $Z_t^F=\cup_{i\in F}(\eta^i_t)$ and
$Z_t^G=\cup_{i\in G}(\eta^i_t)$. Denote by $\widetilde{\eta}_0^i$
the residual of $\eta_0^i$ with respect to $H$ and by
$f$ and $g$ the cardinalities respectively of $F$ and $G$.
Then, by (\ref{assurdo}), we obtain
\begin{equation}\label{seconda}
h_{\P^{\mathbf{n}}}(\Phi^2\cup\Sigma^2\cup Z^F_0\cup
Z^G_t,\mathbf{a})
<\left(1+\sum_{i=1}^k n_i\right)(s-\varepsilon)+2\varepsilon.
\end{equation}

On the other hand, by the semicontinuity of the Hilbert function
there exists an open neighborhood $O$ of $0$ such that for any $t\in
O$
\[
h_{\P^{\mathbf{n}}}\left(\Phi\cup\Sigma^2\cup \left(\cup_{i\in F}\widetilde{\eta}^i_0\right)
\cup Z^G_t,\mathbf{a'}\right)
 \geq   h_{\P^{\mathbf{n}}}\left(\Phi\cup\Sigma^2\cup\left(\cup_{i\in F}\widetilde{\eta}^i_0\right)\cup Z^G_0,\mathbf{a'}\right).
\]
Since $\Phi\cup\Sigma^2\cup\left(\cup_{i\in F}\widetilde{\eta}^i_0\right)\cup Z^G_0
\subseteq \Phi\cup\Sigma^2\cup\Gamma^2_{|H}$,
by Step 2 we compute
\[
h_{\P^{\mathbf{n}}}\left(\Phi\cup\Sigma^2\cup\left(\cup_{i\in F}\widetilde{\eta}^i_0\right)\cup Z^G_0,\mathbf{a'}\right)=s'+\left(1+\sum_{i=1}^k n_i\right)(s-s'-\varepsilon)+f+2g.\]
Since $\Phi^2_{|H}\cup
\left(\cup_{i\in F}\gamma^i\right)$ is a subscheme of $\Phi^2_{|H}\cup\Gamma$, from Step 3
it follows that
\[
h_{\P^{\mathbf{n'}}}\left(\Phi^2_{|H}\cup
\left(\cup_{i\in F}\gamma^i\right),\mathbf{a}\right)\geq s'\sum_{i=1}^k n_i+f
\]
Hence for any $0\neq t\in O$, applying the Castelnuovo inequality
to the scheme $\Omega = \Phi^{2}\cup\Sigma^{2}\cup Z^F_0\cup
Z^G_t$, we get
\begin{eqnarray*}
h_{\P^{\mathbf{n}}}(\Omega,\mathbf{a})
& \geq & h_{\P^{\mathbf{n}}}\left(\Phi\cup\Sigma^2\cup\left(\cup_{i\in F}\widetilde{\eta}^i_0\right)
\cup Z^G_t,\mathbf{a'} \right)+
h_{\P^{\mathbf{n'}}}\left(\Phi^2_{|H}\cup\left(\cup_{i\in F}\gamma^i\right),\mathbf{a}\right) \\
& \geq & s'+\left(1+\sum_{i=1}^k n_i\right)(s-s'-\varepsilon)+f+2g+s'\sum_{i=1}^k n_i+f  \\
& = & \left(1+\sum_{i=1}^k n_i\right)(s-\varepsilon)+2\varepsilon,
\end{eqnarray*}
which contradicts Inequality (\ref{seconda}).
Thus we completed the proof of the claim.
\end{proof}
\begin{example}
Let $\mathbf{n}=(2,2)$ and let $\mathbf{a}=(4,4)$.
In Example~\ref{counter}, we showed that it is impossible to apply
Theorem~\ref{th:horace-easy} to prove the truth of
$T(\mathbf{n};\mathbf{a};45)$.  In this example, we illustrate how  to
reduce $T(\mathbf{n};\mathbf{a};45)$ to computing the dimensions of
secant varieties of ``smaller" Segre-Veronese varieties using
Theorem~\ref{th:horace}.

Let $s'$ and $\varepsilon$ be the quotient and remainder when dividing
$45(2+2+1)-{5\choose 2}{6 \choose 2}$ by $2+2$ respectively. Then
$s'=18$ and $\varepsilon=3$. Thus $s'$ and $\varepsilon$ clearly
satisfy  $s'>\varepsilon$.
Since
\[
120 = (45-18-3)(2+2+1)>{4 \choose 2}{6 \choose 4}=90,
\]
the $5$-tuple $(2,2;2,4;45-18-3)$ is superabundant.
Thus, by Theorem~\ref{th:horace}, one can reduce $T(2,2;4,4;45)$ to
$T(1,2;4,4;18)$, $T(2,2;3,4; 27)$ and $T(2,2;2,4;24)$.

In order to complete the proof of the truth of $T(2,2;4,4;45)$, one can apply Theorem~\ref{th:horace} to $T(1,2;4,4;18)$, $T(2,2;3,4; 27)$ and $T(2,2;2,4;24)$. Like $T(2,2;4,4;45)$,  each statement will be reduced to three sub-statements, each of which can be reduced to other three sub-statements by applying Theorem~\ref{th:horace}. One must repeat this process until one achieves either the statements that are all known to be true or the statements that are small enough, so that one can computationally prove that they are true as indicated in Remark~\ref{rem:computation}.
For example, we checked the truth of $T(1,2;4,4;18)$, $T(2,2;3,4; 27)$ and $T(2,2;2,4;24)$ directly using {\tt Macaulay2}. This shows the truth of $T(2,2;4,4;45)$.
\end{example}
Unfortunately, if $k=2$ and  if one of $a_i$'s is $1$, then  it is
often impossible to apply Theorem~\ref{th:horace}.  For example, if
$(\mathbf{n};\mathbf{a};s)=(2,2;1,4;9)$, then $s'=3$ and $\varepsilon
= 3$. Thus $15 = (9-3-3)(2+2+1) < {2+2 \choose 2}(2+1)=18$, and so
$(2,2;1,2;3)$ is not superabundant. Therefore, we cannot reduce
$T(2,2;1,4;9)$ to $T(1,2;1,4;3)$, $T(2,2;1,3;6)$ and $T(2,2;1,2;3)$.
Another goal of this section is to provide a different approach to
give a way around this kind of problem.
In Example \ref{segre1} we will explain how to apply this second approach to prove the truth of $T(2,2;1,4;9)$.
\begin{definition}
Let  $a_1 =1$ and let $\pi:\P^{\mathbf{n}} \rightarrow \prod_{i=2}^k
\P^{n_i}$ be the canonical projection. For each point $p \in
\P^{\mathbf{n}}$, let $f_p$ be the double point $p^2$ restricted to
$\pi^{-1}(\pi(p))$.
Consider general points $p_1, \dots, p_s$, $q_1, \dots, q_t$, $r_1, \dots, r_v \in  \P^{\mathbf{n}}$
and let $Z=\{p_1^2, \dots, p_s^2, q_1, \dots, q_t, f_{r_1}, \dots, f_{r_v}\}$.
We say that the statement $S(\mathbf{n};\mathbf{a};s;t;v)$ is true if
$T(\mathbf{n};\mathbf{a};Z)$ is true, that is, if
\[
h_{\P^{\mathbf{n}}}(Z, \mathbf{a}) =
 \min \left\{
 s\left(1+\sum_{i=1}^k n_i\right)+t+v(n_1+1), \ N
 \right\}.
\]
We will also write $(\mathbf{n};\mathbf{a};s;t;v)$ for $(\mathbf{n}; \mathbf{a};Z)$.
\end{definition}
\begin{remark}
\label{rem:remark}
Let $\mathbf{n}$ and $\mathbf{a}$ be $k$-tuples of non-negative
integers. We make the following simple remarks:
\begin{itemize}
\item[(i)] $S(\mathbf{n};\mathbf{a};s;0;0)$ is true if and only if $T(\mathbf{n};\mathbf{a};s)$ is true.
\item[(ii)] If $(\mathbf{n};\mathbf{a};s;t;v)$ is subabundant and if
  $S(\mathbf{n};\mathbf{a};s;t;v)$ is true, then
  $(\mathbf{n};\mathbf{a};s';t';v')$ is subabundant and
  $S(\mathbf{n};\mathbf{a};s';t';v')$ is true for any choice of $s'$,
  $t'$ and $v'$ with
$s' \leq s$, $t' \leq t$ and $v' \leq v$.
\item[(iii)] If $(\mathbf{n};\mathbf{a};s;t;v)$ is superabundant and
  if the statement $S(\mathbf{n};\mathbf{a};s;t;v)$ is true, then
  $(\mathbf{n};\mathbf{a};s';t';v')$ is superabundant  and
  $S(\mathbf{n};\mathbf{a};s';t';v')$ is true for any choice of $s'$,
  $t'$ and $v'$ with
$s \leq s'$, $t \leq  t'$ and $v \leq v'$.
This implies that if  $\underline{s} = \left\lfloor \frac{\prod_{i=1}^k {n_i+a_i \choose a_i}}{1+\sum_{i=1}^kn_i} \right\rfloor$
and
$\overline{s} =
\left\lceil \frac{\prod_{i=1}^k {n_i+a_i \choose a_i}}{1+\sum_{i=1}^kn_i} \right\rceil$,
then, in order to prove the truth of
$T(\mathbf{n};\mathbf{a})$,  it is sufficient to show that
$T(\mathbf{n};\mathbf{a};s)$ are true for both $s \in \{\underline{s},\overline{s}\}$.
\item[(iv)]
The following statements are equivalent and  have the same abundancy:
\begin{itemize}
 \item[-] $S(0,\mathbf{n};1,\mathbf{a};s;t;v)$.
 \item[-]  $S(0,\mathbf{n};1,\mathbf{a};s;t+v;0)$.
\item[-]  $S(\mathbf{n};\mathbf{a};s;t+v;0)$.
\end{itemize}
\item[(v)]
If $(\mathbf{n};\mathbf{a};s;t;0) $ is subabundant, then it is
clear that, since the $t$ simple points are assumed to be general,
$S (\mathbf{n};\mathbf{a};s;0;0)=T(\mathbf{n};\mathbf{a};s)$ is
true if and only if $S(\mathbf{n};\mathbf{a};s;t;0)$ is  true.
\end{itemize}
\end{remark}
The following theorem describes the induction procedure we can apply to
study Segre-Veronese varieties when one of the degree is one. This
technique is inspired by the paper \cite{AOP2}, where the authors
study Segre varieties.
\begin{theorem}
\label{th:segrein}
Let $a_1 =1$, $n_1=n_1'+n_1''+1$, $s=s'+s''$ and $t=t'+t''$,  and let
$\mathbf{n}'=(n_1',n_2, \dots, n_k)$, $\mathbf{n}''=(n_1'',n_2, \dots, n_k) \in \N^k$.
Suppose that $(\mathbf{n}';\mathbf{a};s';t';v+s'')$ and
$(\mathbf{n}'';\mathbf{a};s'';t'';v+s')$ are subabundant
(resp. superabundant). 
If $S(\mathbf{n}';\mathbf{a};s';t';v+s'')$ and $S(\mathbf{n}'';\mathbf{a};s'';t'';v+s')$ are true,
then $(\mathbf{n};\mathbf{a};s;t;v)$ is subabundant
(resp. superabundant) 
and $S(\mathbf{n};\mathbf{a};s;t;v)$  is true.
\end{theorem}
\begin{proof}
We only focus on the case when
$(\mathbf{n}';\mathbf{a};s';t';v+s'')$ and
$(\mathbf{n}'';\mathbf{a};s'';t'';v+s')$ are subabundant, because the
remaining case can be proved in a similar fashion.

Let $U$ be a $(n_1'+1)$-dimensional subspace of $V_1$.
Then we have the following Koszul complex:
\[
 \cdots
  \rightarrow (V_1/U)^* \otimes  \ko_{\P^{\mathbf{n}}}(\mathbf{a}') \rightarrow \ko_{\P^{\mathbf{n}}}(\mathbf{a})
 \rightarrow \ko_{\P^{\mathbf{n}'}}(\mathbf{a}) \rightarrow 0,
\]
where
$\mathbf{a}=(1, a_2, \dots, a_k)$.
Let $\iota : U \rightarrow V_1$ be the inclusion. The linear transformation from
$H^0( \ko_{\P^{\mathbf{n}}}(\mathbf{a})) = V_1^* \otimes H^0(\ko_{\P^{\mathbf{n}}}(\mathbf{a}'))$  to $H^0(\ko_{\P^{\mathbf{n}'}}(\mathbf{a}))  = U^* \otimes H^0(\ko_{\P^{\mathbf{n}}}(\mathbf{a}'))$ induced by the last map of the Koszul complex is given by $\iota^* \otimes \mathrm{id_{H^0(\ko_{\P^{\mathbf{n}}}(\mathbf{a}'))}}$, and hence it is surjective.
By taking the cohomology, we therefore obtain the following short exact sequence:
\[
 0 \rightarrow (V_1/U)^* \otimes H^0(\ko_{\P^{\mathbf{n}}}(\mathbf{a}'))
 \rightarrow   H^0( \ko_{\P^{\mathbf{n}}}(\mathbf{a}))
 \rightarrow  H^0(\ko_{\P^{\mathbf{n}'}}(\mathbf{a}))
\rightarrow  0.
\]
Taking the dual of the first linear transformation of the above sequence
yields the rational map $\varphi$ from  $\prod_{i=1}^k \P^{n_i}$ to
$\P^{\mathbf{n}''}=\P(V_1/U) \times \prod_{i=2}^k \P^{n_i}$.

Let $Z=\{p_1^2, \dots, p_s^2\}$, let $\Phi=\{q_{1}, \dots, q_{t}\}$ and let $\Psi=\{f_{r_1}, \dots, f_{r_v}\}$.
Suppose that $\{p_1, \dots, p_{s''} \}$ and $\{q_1, \dots, q_{t''}\}$
are not contained in  $\P^{\mathbf{n}'}$, but the rest of the $p_i$'s
and  $q_i$'s are in $\P^{\mathbf{n}'}$, while the $r_i$'s are general points.
Then we have the following short exact sequence:
\[
 0 \rightarrow \sI_{Z \cup \Phi \cup \Psi \cup \P^{\mathbf{n}'}}(\mathbf{a})
 \rightarrow \sI_{Z \cup \Phi\cup \Psi}(\mathbf{a}) \rightarrow
 \sI_{(Z \cup \Phi \cup \Psi) \cap \P^{\mathbf{n}'}, \P^{\mathbf{n}'}}(\mathbf{a})
 \rightarrow 0.
\]
Let $\pi''$ be the canonical projection from $\P^{\mathbf{n}''}$ to $\prod_{i=2}^k \P^{n_i}$
and let $Z''$ be the following zero-dimensional subscheme of $\P^{\mathbf{n}"}$:
\[
\{\varphi(p_1)^2, \dots, \varphi(p_{s''})^2, \varphi(q_1), \dots,
\varphi(q_{t''}), f_{\varphi(r_1)}, \dots, f_{\varphi(r_{v})},
f_{\varphi(p_{s''+1})},\ldots, f_{\varphi(p_{s})} \}.
\]
One can immediately show that $H^0(\P^{\mathbf{n}''}, \sI_{Z''}(\mathbf{a}))$ is isomorphic to $H^0\left(\P^{\mathbf{n}}, \sI_{Z \cup \Phi \cup \P^{\mathbf{n}'}}(\mathbf{a})\right)$.

Let $\psi$ be the projection from $\P^{\mathbf{n}} \setminus \P^{\mathbf{n}''}$ to $\P^{\mathbf{n}'}$ and
let $Z'$ be the following zero-dimensional subscheme  of $\P^{\mathbf{n}'}$:
\[
\{p_{s''+1}^2,\dots,p_s^2, q_{t''+1}, \dots, q_t, f_{\psi(r_{1})}, \dots, f_{\psi(r_v)}, f_{\psi(p_1)}, \dots, f_{\psi(p_{s''})}\}.
\]
Note that  $H^0(\sI_{Z'}(\mathbf{a}))$ is isomorphic to
$\left(I_{Z
  \cup
  \Phi}+I_{\P^{\mathbf{n}'}}/I_{\P^{\mathbf{n}'}}\right)_{\mathbf{a}}$.
This implies that if
\[
 h_{\P^{\mathbf{n}''}}(Z'',\mathbf{a})= s''\left(1+n_1''+\sum_{i=2}^k n_i \right)+t''+(v+s')(n_1''+1)
\]
and
\[
h_{\P^{\mathbf{n}'}}(Z', \mathbf{a})
= s'\left(1+n_1'+\sum_{i=2}^kn_i\right)+t'+(v+s'')(n_1'+1) \]
then $h_{\P^{\mathbf{n}}}\left(Z \cup \Phi,\mathbf{a}
\right)=s\left(1+\sum_{i=1}^k n_i\right)+t+v(n_1+1)$, which completes
the proof.
\end{proof}
\begin{example}\label{segre1}
As the first application of Theorem~\ref{th:segrein}, we will show
that $T(\mathbf{n};\mathbf{a};s)$ is true with
$(\mathbf{n};\mathbf{a};s)=(2,2;1,4;9)$.  Note that
$(\mathbf{n};\mathbf{a};s)$ is subabundant.
Let $s'=6$. Then $s''=9-6=3$. Since $(0,2;1,4;3;0;6)$ and
$(1,2;1,4;6;0;3)$ are equiabundant, we can reduce
$T(\mathbf{n};\mathbf{a};s)$ to $S(0,2;1,4;3;0;6)=S(2;4;3;0;6)$ and
$S(1,2;1,4;6;0;3)$. The statement $S(1,2;1,4;6;0;3)$ can be reduced to
twice \linebreak
$S(0,2;1,4;3;0;6)$.
In order to prove that $T(\mathbf{n};\mathbf{a};s)$ is true, it is
therefore enough to prove the truth of $S(0,2;1,4;3;0;6)$. Note that
$S(0,2;1,4;3;0;6)$ and  $S(2;4;3;6;0)$ are the same statements by Remark \ref{rem:remark} (iv). Also,
the condition that $S(2;4;3;6;0)$ is true is equivalent to the
condition that $S(2;4;3;0;0)=T(2;4;3)$ is true
by Remark \ref{rem:remark} (v). It is known by the
Alexander-Hirschowitz theorem that $T(2;4;3)$ is true. Thus
$T(\mathbf{n};\mathbf{a};s)$ is also true.
\end{example}
Let $\mathbf{n}$, $\mathbf{a} \in \N^k$. As already stated in
Section~\ref{sec:intro},  Theorem~\ref{th:horace} cannot be applied to
any secant variety of $X_{\mathbf{n}, \mathbf{a}}$ if
$\mathbf{a}=(2^k)$.  Theorem~\ref{th:segrein} cannot be used directly
in this case either. In the following example, we illustrate how to
combine an argument based on the Castelnuovo inequality with
Theorem~\ref{th:segrein} to study secant varieties of such
Segre-Veronese varieties:
\begin{example}
\label{segre2}
Here we prove that $T(2,2;2,2;5)$ is true.  Let $p_1, \dots, p_5$ be
generic points of $(\P^2)^2$ and let $Z = \{p_1^2, \dots,
p_5^2\}$. Specializing $p_1$, $p_2$ and $p_3$ to $H = \P^1 \times \P^2
\subset (\P^2)^2$ yields a short exact sequence
\[
 0 \rightarrow \sI_{\widetilde{Z}}(1,2) \rightarrow
 \sI_Z(2,2) \rightarrow \sI_{Z \cap H, H}(2,2)
 \rightarrow 0.
\]
It was shown by Bauer and Draisma~\cite{BD} that $h_{H}(Z \cap H,
(2,2))$ has the expected value, i.e., $T(1,2;2,2;3)$ is true. It
suffices therefore to show that $\widetilde{Z}$ has the expected value at $(1,2)$.
Note that $\widetilde{Z}=\{p_1, p_2, p_3, p_4^2, p_5^2\}$. Recall that
$p_1$, $p_2$ and $p_3$ lie in $H$. Thus specializing $p_5$ to $H$,  we
can reduce the above-mentioned statement to $S(1,2;1,2;1;3;1)$ and
$S(0,2;1,2;1;0;1)$. Note that $S(0,2;1,2;1;0;1)$ is equivalent to
$S(2;2;1;1;1;0)$.  Since $S(2;2;1;1;0;0)$ is true, so is
$S(2;2;1;1;1;0)$ by Remark~\ref{rem:remark}. Thus it remains to show
that $S(1,2;1,2;1;3;1)$ is true. This statement  can be reduced to
$S(0,2;1,2;1;1;1)$ and $S(0,2;1,2;0;2;2)$.
By Remark~\ref{rem:remark}, $S(0,2;1,2;1;1;1)$ and $S(0,2;1,2;0;2;2)$
are equivalent to $S(2;2;1;2;0)$ and $S(2;2;0;4;0)$ respectively.
Clearly, the latter statement is true. Also, since \linebreak
$S(2;2;1;0;0)$ is
true, so is $S(2;2;1;2;0)$.
Thus $S(1,2;1,2;1;3;1)$ is true.
Therefore, $T(1,2;2,2;3)$ is true.
\end{example}
We conclude this section by presenting immediate, but useful
consequences of Theorem~\ref{th:horace-easy} and Lemma~\ref{chandler}.
\begin{lemma}\label{small-s}
Let $\mathbf{a},\mathbf{b},\mathbf{n},\mathbf{m}\in (\Z_{\geq 0})^k \setminus \{(0,\dots, 0)\}$ and let $s\in\N$.
Suppose that  $\mathbf{a}\leq\mathbf{b}$ and $\mathbf{n}\leq\mathbf{m}$.
\begin{itemize}
\item[(i)] If $T(\mathbf{n}; \mathbf{a}; s)$ is true, if $(\mathbf{n}; \mathbf{a}; s)$ is subabundant and
if  $\mathbf{a} \geq (1, \dots, 1)$,
then $T(\mathbf{n}; \mathbf{b}; s)$ is true and $(\mathbf{n}; \mathbf{b}; s)$ is subabundant.
\item[(ii)] If $T(\mathbf{n}; \mathbf{a}; s)$ is true, if $(\mathbf{n}; \mathbf{a}; s)$ is subabundant
and if
\[
s\le \binom{n_\ell+a_\ell-1}{a_\ell-1}\prod_{i\neq \ell} {n_i+a_i \choose a_i}
\]
for all $\ell$ such that $m_\ell>n_\ell$ and $a_{\ell} \geq 1$,
then $T(\mathbf{m}; \mathbf{a}; s)$ is true and $(\mathbf{m}; \mathbf{a}; s)$ is subabundant.
\item[(iii)] If $T(\mathbf{n}; \mathbf{b}; s)$ is true and if $(\mathbf{n}; \mathbf{b}; s)$ is superabundant,
then $T(\mathbf{n}; \mathbf{a}; s)$ is true and $(\mathbf{n}; \mathbf{a}; s)$ is superabundant.
\end{itemize}
\end{lemma}
\begin{proof}
Note that if $(\mathbf{n}; \mathbf{a}; s)$ is subabundant, then so are $(\mathbf{n}; \mathbf{b}; s)$ and $(\mathbf{m}; \mathbf{a}; s)$.

\vspace{1mm}
\noindent
(i) Since $\mathbf{a} \geq (1, \dots, 1)$, without loss of generality, we may assume that
$1\leq a_1 < b_1$.
Let $\mathbf{a'}=(a_1+1,a_2,\ldots,a_k)$ and $H$
a hyperplane defined by a linear form in $R_{(1,0,\ldots,0)}$.
By induction it suffices to prove that $T(\mathbf{n};\mathbf{a'};s)$ is true.
Consider a collection $Z$ of $s$ general double points in $\P^{\mathbf{n}}$.
Suppose that the support of $Z$ is not contained in $H$.
>From the short exact sequence,
\[
0\rightarrow \sI_{Z}(\mathbf{a})\rightarrow
\sI_Z(\mathbf{a'})\rightarrow \mathcal{O}_H(\mathbf{a'})\rightarrow 0,
\]
we can conclude that
$h_{\P^{\mathbf{n}}}(Z,\mathbf{a'})$ is the expected value, because the trace of $Z$ is empty and $\widetilde{Z} = Z$.

\vspace{1mm}
\noindent
(ii) The statement is trivial if $\mathbf{n} = \mathbf{m}$. Thus we
may assume that $\mathbf{n} < \mathbf{m}$. Then there exists at least
one $\ell \in \{1, \dots, k\}$ such that $m_{\ell} > n_{\ell}$ and
$a_{\ell} \geq 1$, because otherwise $T(\mathbf{n};\mathbf{a};s)$ and
$T(\mathbf{m};\mathbf{a};s)$ are the same statement.
Without loss of generality, we may assume that $\ell=1$.
Then, by induction, it is enough to prove that
$T(\mathbf{n'};\mathbf{a};s)$ is true for $\mathbf{n'}=(n_1+1,n_2,\ldots,n_k)$.
Consider a collection $Z$ of $s$ general double points of $\P^{\mathbf{n'}}$.
Suppose that they are all  contained in
$H=\P^{\mathbf{n}}\subset\P^{\mathbf{n'}}$.
Hence the trace of $Z$ is given by $s$ double
points of $H$, while the residual $\widetilde{Z}$
 is given by $s$ simple points
contained in the hyperplane $H$.
Then we have the following exact sequence
\[
0\rightarrow \sI_{\widetilde{Z}}(a_1-1,a_2,\ldots,a_k)\rightarrow
\sI_Z(\mathbf{a})\rightarrow \sI_{Z\cap H}(\mathbf{a})\rightarrow 0.
\]
By assumption, $T(\mathbf{n}; \mathbf{a}; s)$ is true. Thus (ii)
immediately follows from Lemma \ref{chandler} and from the assumption that
$s\le \binom{n_1+a_1-1}{a_1-1}\prod_{i=2}^k {n_i+a_i \choose a_i}$.

\vspace{1mm}
\noindent
(iii)
Clearly if $(\mathbf{n}; \mathbf{b}; s)$ is superabundant,
then $(\mathbf{n}; \mathbf{a}; s)$ is also superabundant.
Given $\mathbf{b'}=(b_1-1,b_2,\ldots,b_h)$, we only need to prove $T(\mathbf{n};\mathbf{b'};s)$.
As in the proof of (i), we consider a collection $Z$ of $s$ general
double points whose support is not contained in $H$.  Then, by the
Castelnuovo exact sequence,  we can immediately see that
$T(\mathbf{n};\mathbf{b'};s)$ is true.
\end{proof}
\section{Two-factor Segre-Veronese varieties}
\label{sec:two-factor}
The purpose of this section is to establish the existence of a class
of non-defective two-factor Segre-Veronese varieties. First of all, we
will recall  some basic results on  secant varieties of such
Segre-Veronese varieties.  Let $\mathbf{n}=(m,n)$, $\mathbf{a}=(a,b)
\in \mathbb{N}^2$ and let $N (\mathbf{n},\mathbf{a})= {m+a \choose
  a}{n+b \choose b}$. We use just $N$ instead of
$N(\mathbf{n};\mathbf{a})$ if $\mathbf{n}$ and $\mathbf{a}$ are clear
from the context. As in the previous section, we denote by
$X_{\mathbf{n},\mathbf{a}}$ the Segre-Veronese variety obtained from
$\P^\mathbf{n}$ by embedding in $\P^{N-1}$ by the morphism given by
$\cO(\mathbf{a})$.  Let $\underline{s}(\mathbf{n},\mathbf{a}) =
\lfloor N/(m+n+1) \rfloor$ and let
$\overline{s}(\mathbf{n},\mathbf{a}) = \lceil N/(m+n+1) \rceil$.  We
write $\underline{s}$ and $\overline{s}$ instead of
$\underline{s}(\mathbf{n},\mathbf{a})$ and
$\overline{s}(\mathbf{n},\mathbf{a})$ respectively if $\mathbf{n}$ and
$\mathbf{a}$ are clear from the context.
As mentioned in Remark~\ref{rem:remark}, in order to prove
that $T(\mathbf{n};\mathbf{a})$ is true,  it is sufficient to show
that $T(\mathbf{n};\mathbf{a};s)$ for $s = \underline{s}$ and
$\overline{s}$.

As was mentioned earlier, the problem of finding the dimension of
$\sigma_s(X_{\mathbf{n},\mathbf{a}})$ can be translated into the
problem of calculating the value of the multi-graded Hilbert function
of $s$ double points on $\P^{\mathbf{n}}$ at $\mathbf{a}$.
In their several papers,  Catalisano, Geramita, and Gimigliano showed
the relationship between ideals of varieties in multi-projective space
and ideals in standard polynomial rings.
In \cite[Theorem 2.1]{CGG2}, they used it to prove the following theorem:
\begin{theorem}[\cite{CGG2}]
\label{th:CGG}
$T(1,1;a,b;s)$ is true except for $a=2$, $b=2d$ ($d \geq 1$) and $s=b+1$.
\end{theorem}
This theorem was also proved by Baur and Draisma. Their proof uses
tropical techniques (see \cite[Theorem 1.1]{BD} for more details).
\begin{example}\label{example-easy}
As the first  application of our techniques we prove that
$T(m,1;1,2;2)$ is true for any $m\ge1$.

By Theorem~\ref{th:CGG},  $T(1,1;1,2;2)$ is true. Moreover $(1,1;1,2;2)$ is
equiabundant. Since
\[
s=2< 3= \binom{n_1+a_1-1}{a_1-1} {n_2+a_2 \choose a_2},
\]
we can deduce that $T(m,1;1,2;2)$ are true for all $m\ge1$ by Lemma \ref{small-s} (ii).
\end{example}
Let $s$ be a positive integer and let $s'$ and $\varepsilon$ be the
quotient and remainder when dividing $s(m+n+1)-{m+a-1 \choose a-1}
{n+b \choose b}$ by $m+n$. In order to prove
the truth of $T(m,n;a,b;s)$, we need to show that
the $5$-tuple  $(m,n;a-2,b;s-s'-\varepsilon)$ is superabundant.
The following lemma proves that this is actually the case  for most of $(m,n;a,b)$.
\begin{lemma}
Let $a, b \geq 3$. For each $1 \leq s \leq \left\lceil\frac{{m+a
    \choose a}{n+b \choose n}}{m+n+1}\right\rceil$, let $s'$ and
$\varepsilon$ be as above. Then $(m,n;a-2,b;s-s'-\varepsilon)$ is
superabundant unless $(m,n)=(1,1)$. 
\label{th:supab}
\end{lemma}
\begin{proof}
We want to prove that the integer $F(m,n;a,b)$ is non-negative, where
\begin{eqnarray}
 F(m,n;a,b)=(s-s'-\varepsilon)(m+n+1)-N(m,n;a-2,b).
\label{seq:superabundant}
\end{eqnarray}
By definition,
\begin{equation}\label{definitions'}
s\,(m+n+1)-N(m,n;a-1,b) =s'\,(m+n)+\varepsilon,
\end{equation}
where $0  \leq \varepsilon \leq  m+n-1$.
So we have
\begin{eqnarray*}
&& F(m,n;a,b) \\
&=& s\,(m+n+1)-s'\,(m+n)-\varepsilon-s'-\varepsilon\,(m+n)-N(m,n;a-2,b) \\
&=& N(m,n;a-1,b)-N(m,n;a-2,b)-s'-\varepsilon\,(m+n) \\
&=& N(m-1,n;a-1,b)-s'-\varepsilon\,(m+n)
\end{eqnarray*}
Since $s\,(m+n+1) \leq N(m,n;a,b)+(m+n+1)$ by assumption,  the
following inequality holds:
\[
s\,(m+n+1)- N(m,n;a-1,b) \leq N(m-1,n;a,b) +(m+n+1).
\]
This implies that $s'\,(m+n) \leq N(m-1,n;a,b)+(m+n+1)$, i.e.,
\[
s' \leq \frac{1}{m+n}\left\{N(m-1,n;a,b)+(m+n+1)\right\}.
\]
Thus we obtain
\begin{eqnarray*}
& & F(m,n;a,b) \\
 & \geq & N(m-1,n;a-1,b)-
\frac{\left\{N(m-1,n;a,b)+(m+n+1)\right\}}{m+n}
-\varepsilon\,(m+n) \\
& \geq & \frac{{n+b \choose n}}{m+n}H(m,n;a)-\frac{m+n+1}{m+n}-(m+n-1)(m+n),
\end{eqnarray*}
where  $H(m,n;a) =
\frac{(m+a-2)!}{a!(m-1)!}\left\{a(m+n)-(m+a-1)\right\}$. Note that
$H(m,n;a)$ is an increasing function of $a$ if $a \geq 3$.   Let
\[
G(m,n;a,b)=\frac{{n+b \choose n}}{m+n}H(m,n;a)-\frac{m+n+1}{m+n}-(m+n-1)(m+n).
\]
It follows that $G(m,n;a,b)$ is an increasing function of $a$ and $b$, if $a\ge3$.
It is not very hard to show that $G(m,n;3,3) \geq -\frac{2}{3}$
unless $(m,n)=(1,1)$.
Hence we have $F(m,n;a,b)\geq G(m,n;a,b)\ge G(m,n;3,3)\geq -\frac23$.
Since $F(m,n;a,b)$ is an integer, we can conclude that is it non-negative.
\end{proof}
In the following lemma, we show that the inequality $s'\ge\varepsilon$
holds in most cases:
\begin{lemma}\label{s'big}
Let $a\ge3$ and let  $b,m,n\ge1$.
For each $s\ge \left\lfloor\frac{{m+a\choose{a}}{n+b\choose
      n}}{m+n+1}\right\rfloor$, let $s'$ and
$\varepsilon$ be as above. Then $s'\ge \varepsilon$
in the following cases:
\begin{itemize}
\item[(i)] $b\ge3$;
\item[(ii)] $b=1$ and $m\ge3$;
\item[(iii)] $b=1$, $m=2$ and $n=1$.
\end{itemize}
\end{lemma}
\begin{proof}
Since $n+m-1\ge \varepsilon$, it suffices to show that $s'\ge n+m-1$.
Assume that $s' < n+m-1$.
By assumption we know that
$s(m+n+1)\ge {m+a\choose a}{n+b \choose n}-(m+n)$. Combining this
relation with \eqref{definitions'} yields
\[
N(m,n;a,b)-(m+n)\le s'\,(m+n)+\varepsilon+N(m,n;a-1,b)
\]
from which we obtain
\begin{eqnarray*}
N(m-1,n;a,b)
& \le &
s'\,(m+n)+\varepsilon+(m+n) \\
& \le &
(n+m-2)(m+n)+(m+n-1)+(m+n) \\
& = &
(n+m)^2-1.
\end{eqnarray*}
Now we need to prove that this inequality provides a contradiction in each case.
Let $G(m,n,a,b)=N(m-1,n;a,b)-(m+n)^2$. It is enough to prove that $G(m,n,a,b) > 0$.

\vspace{2mm}
\noindent
(i) Suppose that $b  \geq 3$.
Note that
$N(m-1,n;a,b) \geq N(m-1,n;3,3)$.
It follows therefore that  if $a,b\ge3$, then $G(m,n,a,b)\ge G(m,n;3,3)$.
It is straightforward to prove that $G(m,n,3,3)$ is positive for all $m,n\ge1$.

\vspace{2mm}
\noindent
(ii) Suppose that $b = 1$ and $m \geq 3$.
In the same way as in (i), one can prove that
$G(m,n,a,b)\ge G(m,n;3,1)$. It is not hard to show that
$G(m,n,3,1)$ is positive when $m\ge3,$ and $n\ge1$.

\vspace{2mm}
\noindent
(iii)  Assume that  $b=1,m=2$ and $n=1$. Then we have
$\varepsilon\le 2$. We want to prove that $s'\ge 2$. Assume for the contradiction
that $s'\le1$. By the hypothesis we have $s\ge
\left\lfloor\frac{(a+2)(a+1)}{4}\right\rfloor$, which implies
$4s\ge (a+2)(a+1)-3$. By \eqref{definitions'} we have
\[
(a+2)(a+1)-3\le a(a+1)+3s'+\varepsilon\le a(a+1)+3+2,
\]
or  $2(a-3)\le0$, which is false for all $a\ge4$. If $a=3$,
then we have $s\ge5$. On the other hand,  \eqref{definitions'} gives rise to
$4s\le 12+3+2=17$, which is a contradiction.
\end{proof}
The result presented below was already proved by Chiantini and Ciliberto \cite{ChiCil}. Here we give a different proof to illustrate how the Horace method works.
\begin{theorem}
\label{n11d}
$T(n,1;1,d)$ is true for any $n,d\ge1$.
\end{theorem}
\begin{proof}
The proof is by induction on $d$. It is immediate to check that $T(n,1;1,1)$ is
true (see Section~\ref{sec:intro}). The truth of the statement $T(n,1;1,2)$ immediately follows
from Example \ref{example-easy} and \cite[Example 2.9]{AB}. Thus
we may assume $d\geq 3$ and $n\geq1$.

We first prove the truth of $T(n,1;1,d;s)$ for
$s=\underline{s}(n,1;1,d)=\left\lfloor\frac{(n+1)(d+1)}{n+2} \right\rfloor$.
Let $p_1, \dots, p_s$ be points on $\P^n \times \P^1$ and let
$Z=\{p_1^2, \dots, p_s^2\}$. Suppose that  $p_s$ lies in a hyperplane
$H$ of degree $(0,1)$. Then we get the following sequence:
\[
 0 \rightarrow \sI_{\tZ}(1,d-1) \rightarrow
 \sI_Z(1,d) \rightarrow \sI_{Z \cap H, H}(1,d)
 \rightarrow  0,
\]
where $\tZ=\{p_1^2, \dots, p_{s-1}^2\}\cup \{p_s\}$.
Since  the trace of $Z$ consists of only one double point of $H$, we have
\[
h_H(Z \cap H, (1,d))=n+1.
\]
By induction hypothesis,  $T(n,1;1,d-1;s-1)$ and $T(n,1;1,d-2;s-1)$
are both true, and thus
\[
h_{\P^{\mathbf{n}}}(\{p_1^2, \dots, {p_{s-1}}^2\},(1,d-1))=\min\{(s-1)(n+2),(n+1)d\}=(s-1)(n+2)
\]
and
\[
h_{\P^{\mathbf{n}}}(\{p_1^2, \dots,{p_{s-1}}^2\},(1,d-2))=\min\{(s-1)(n+2),(n+1)(d-1)\}.
\]
It is straightforward to prove the inequality
\[
(s-1)(n+2)+1\le\min \{(s-1) (n+2) ,(n+1)(d-1)\}+(n+1).
\]
So it follows from Lemma~\ref{chandler}  that
$h_{\P^{\mathbf{n}}}(\tZ,(1,d-1))= (s-1)(n+2)+1$.
By Theorem \ref{th:horace-easy} we can deduce that $h_{\P^{\mathbf{n}}}(Z,(1,d))=s(n+2)$, because
$h_{\P^{\mathbf{n}}}(\tZ,\mathbf{a'})$ and $h_{\P^{\mathbf{n'}}}(Z\cap H,\mathbf{a})$ are the expected values and they are both subabundant. Thus,  $T(n,1;1,d;s)$ is true.

In a similar manner, we can prove that $T(n,1;1,d;s)$ is true for
$s=\overline{s}(n,1;1,d)$. Let $p_1, \dots, p_s$ be points on $\P^n\times
\P^1$ and let $Z=\{p_1^2, \dots, p_s^2\}$.
Specializing $p_s$ to $H$ yields
the following sequence:
\[
 0 \rightarrow \sI_{\tZ}(1,d-1) \rightarrow
 \sI_Z(1,d) \rightarrow \sI_{Z\cap H, H}(1,d)
 \rightarrow  0.
\]
As in the previous case,  we have $h_H(Z \cap H, (1,d))=n+1$. By induction hypothesis,
$T(n,1;1,d-1;s-1)$ is true. Additionally,
$(n,1;1,d-1;s-1)$ is superabundant. Therefore,
$h_{\P^{\mathbf{n}}}(Z,(1,d))=(n+1)(d+1)$, which completes the proof.
\end{proof}
We recall now a result proved by Abrescia.
\begin{theorem}[\cite{Abr}]\label{abrescia}
$T(n,1;2,2d+1)$  is true for any $n\ge1$ and $d\ge0$.
\end{theorem}
The following is the first application of the differential Horace lemma:
\begin{theorem} \label{th:n1a2d}
$T(n,1;a,2d+1)$ is true for any $d,n,a\ge1$.
\end{theorem}
\begin{proof}
The proof is by double induction on $n$ and $a$.
We know that $T(n,1;1,2d+1)$ is true by Theorem \ref{n11d} and that
$T(n,1;2;2d+1)$ is true by Theorem \ref{abrescia}.
The statement $T(1,1;a,2d+1)$ is also true by Theorem \ref{th:CGG}.

Suppose now that $a\ge 3$ and $n\ge2$.
Recall that it is enough to prove $T(n,1;a,2d+1;s)$ for $s\in\{\underline{s},\overline{s}\}$.
We want to apply Theorem~\ref{th:horace}.
Let $s'$ and $\epsilon$ be the quotient and remainder when dividing
$s(n+2)- N(n,1;a-1,2d+1)$ by $n+1$.
Note that, by Lemma~\ref{th:supab},
$(n,1;a-2,2d+1;s-s'+\varepsilon)$ is superabundant, because $n \geq 2$, $a\ge3$ and $2d+1\ge3$.
Additionally, by Lemma~\ref{s'big}~(i),  we obtain
$s'\ge \varepsilon$.
Now, by induction hypothesis,
$T(n-1,1;a,2d+1)$, $T(n,1;a-1,2d+1)$ and $T(n,1;a-2,2d+1)$
are all true. Thus Theorem~\ref{th:horace} implies that $T(n,1;a,2d+1)$ is true.
\end{proof}
The following theorem is a consequence of Theorem \ref{th:segrein}:
\begin{theorem}
\label{th:n124s} For any $n,d\ge1$, $T(n,1;2,2d;s)$ is true if $s
\leq d(n+1)$ or $s \geq (d+1)(n+1)$.
\end{theorem}
\begin{proof}
To prove this theorem, we only need to show that
$T(n,1;2,2d;d(n+1))$ and $T(n,1;2,2d;(d+1)(n+1))$ are true. The
proof is by induction on $n$. Recall that $T(1,1;2,2d;s)$ is true
unless $s=2d+1$ by Theorem~\ref{th:CGG}.
Also, $T(n,1;1,2d)$ is true by Theorem~\ref{n11d}.

We first prove that  $T(n,1;2,2d;d(n+1))$ is true. Let $s=dn+d$,
$s'=dn$, $s''=d$ and let $H$ be a hyperplane of multi-degree
$(1,0)$.

Specializing $s'$ points to $H$,
since $(n-1,1;2,2d;s')$ and
$(n,1;1,2d;s'';s';0)$ are both subabundant, we can apply
Theorem~\ref{th:horace-easy}.
By induction hypothesis, $T(n-1,1;2,2d;s')$ is true and so
it suffices to prove that
$h_{\P^n\times\P^1}(\widetilde{Z},(1,2d))$ is the expected value,
where $\widetilde{Z}$ is given by
$s''$ general double points and $s'$ simple points contained in $H$
(and general).
In order to prove this fact we apply now Theorem~\ref{th:segrein}.
Since $(n-1,1;1,2d;0;s';s'')$ and $(0,1;1,2d;s'';0;0)$ are both
subabundant, it is enough to prove that
$S(n-1,1;1,2d;0;s';s'')$ and $S(0,1;1,2d;s'';0;0)$ are true.

By Theorem~\ref{n11d},
$S(n-1,1;1,2d;s'';0;0)=T(n-1,1;1,2d;s'')$ is true. This implies
that $S(n-1,1;1,2d;0;0;s'')$ is also true. Additionally, the $s'$
points are in general position in $H$. So $S(n-1,1;1,2d;0;s';s'')$
is true. Since $S(0,1;1,2d;d;0;0)=T(1;2d;d)$ is clearly true, the
theorem follows from Theorem~\ref{th:segrein}.

One can prove that $T(n,1;2,2d;(d+1)(n+1))$ is true
by taking $s'=(d+1)n$ and $s''=d+1$ and by replacing ``subabundant" by
``superabundant" in the previous argument.
\end{proof}
\begin{remark}\label{rem-abrescia}
In \cite{Abr} Abrescia proved Theorem \ref{th:n124s} with different techniques.
Moreover she proved that $\sigma_s(X_{(n,1),  (2,2d)})$ is defective for any $d(n + 1) + 1 \le s \le (d + 1)(n + 1) - 1$.
\end{remark}
\begin{lemma}
\label{th:n134s}
$T(n,1;3,4)$ is true for any $n\ge1$.
\end{lemma}
\begin{proof}
To prove this lemma, it is enough to show that $T(n,1;3,4;s)$ is true for $s \in\{ \underline{s},\overline{s}\}$.
Here we only show that $T(n,1;3,4;s)$ is true for $s = \underline{s}$,
because the remaining case follows the same path.
The proof is by induction on $n$. Note that
$T(1,1;3,4)$ is true by Theorem~\ref{th:CGG}.

Suppose now that $n\ge2$. We also assume by induction that $T(n-1,1;3,4)$ is
true. Let $s'$ and $\varepsilon$ be the
quotient and remainder in the division of $s(n+2)-5{n+2 \choose 2}$ by
$n+1$.
Then, in order to apply Theorem~\ref{th:horace}, it is enough to check that
$T(n,1;2,4;s-s')$ and $T(n,1;1,4;s-s'-\varepsilon)$ are true,
that $(n,1;1,4;s-s'-\varepsilon)$ is superabundant and that $s'\le\varepsilon$.

>From Theorem~\ref{n11d} it follows that $T(n,1;1,4)$ is true.
Moreover $(n,1;1,4;s-s'-\varepsilon)$ is superabundant by
Lemma~\ref{th:supab}, because $n\ge2$ and $s'\ge\varepsilon$ by Lemma~\ref{s'big} (i).
By Theorem~\ref{th:n124s},
$T(n,1;2,4;s-s')$ is true if  $s-s' \leq 2(n+1)$.
Hence our task is to show that the inequality $s'\le\varepsilon$ holds.

It is not hard to prove that the inequality holds for $n=2$, and so we may assume that $n \geq 3$. By the definitions of $s$ and $s'$, we have
\begin{small}
\begin{eqnarray*}
2n+2-s+s' & = & 2n+2-\left\lfloor\frac{5{n+3 \choose 3}}{n+2}\right\rfloor +
\left\lfloor \frac{\left\lfloor\frac{5{n+3 \choose3}}{n+2}
\right\rfloor(n+2)-5{n+2\choose2}}{n+1}\right\rfloor  \\
& = & 2n+2+ \left\lfloor \frac{\left\lfloor\frac{5(n+3)(n+1)}{6}
\right\rfloor-5{n+2\choose2}}{n+1}\right\rfloor \\
& > &  2n+2+  \frac{\left\lfloor\frac{5(n+3)(n+1)}{6}
\right\rfloor-5{n+2\choose2}}{n+1}-1\\
& \ge &
2n+1+  \frac{5(n+3)}{6}-\frac{1}{n+1}
-\frac{5(n+2)}{2} \\
& = & \frac{2n^2-7n-15}{6(n+1)}.
\end{eqnarray*}
\end{small}

\noindent
It is straightforward to show that $f(n) = \frac{2n^2-7n-15}{6(n+1)}$ is an increasing function.
Since $f(3) = -3/4$, we can conclude that $2n+2-s+s' \geq 0$. Thus we completed the proof.
\end{proof}
\begin{lemma}
\label{th:n144s}
$T(n,1;4,4)$ is true for any $n\ge1$.
\end{lemma}
\begin{proof}
In order to prove the truth of this statement, it is enough to show that
$T(n,1;4,4;\underline{s})$ and $T(n,1;4,4;\overline{s})$ are true.
Here we only consider the first case, because the remaining case can be proved in a similar fashion.

We use induction on $n$.
Note that $T(1,1;4,4;s)$ is true by Theorem~\ref{th:CGG}.
It can be also proved directly that $T(2,1;4,4;s)$ is true. So we may assume that $n\ge3$.
Let $s'$ and $\varepsilon$ be the quotient and remainder in the
division of $s(n+2)-5{n+3 \choose 3}$ by $n+1$ respectively.
Since $(n,1;2,4;s-s'-\varepsilon)$ is superabundant by Lemma~\ref{th:supab}
and $s'\ge\varepsilon$ by Lemma \ref{s'big} (i),
the statement $T(n,1;4,4;s)$ can be reduced to $T(n-1,1;4,4;s')$, $T(n,1;3,4;s-s')$
and $T(n,1;2,4;s-s'-\varepsilon)$. By induction hypothesis,
$T(n-1,1;4,4;s')$ is true.
It follows from Lemma~\ref{th:n134s} that $T(n,1;3,4;s-s')$ is
true.  Hence it suffices to prove that the inequality $s-s'-\varepsilon \geq 3n+3$ holds by
Theorem~\ref{th:n124s}.

It is not hard to show that the above inequality holds for $n=2$.
Suppose therefore that $n \geq 3$.  Then
\begin {small}
\begin{eqnarray*}
s-s'-\varepsilon
& = &
\left\lfloor\frac{5{n+4 \choose 4}}{n+2} \right\rfloor
-\left\lfloor \frac{\left\lfloor\frac{5{n+4 \choose4}}{n+2}
\right\rfloor(n+2)-5{n+3\choose3}}{n+1}\right\rfloor -\epsilon \\
& = & -\left\lfloor \frac{\left\lfloor\frac{5{n+4 \choose4}}{n+2}
\right\rfloor-5{n+3\choose3}}{n+1}\right\rfloor -\epsilon \\
& \geq &  \frac{-\left\lfloor\frac{5{n+4 \choose4}}{n+2}
\right\rfloor+5{n+3\choose3}}{n+1} -n \\
& \geq & -\frac{5(n+4)(n+3)}{24}+\frac{5(n+3)(n+2)}{6} - n \\
& = & \frac{15n^2+41n+60}{24}.
\end{eqnarray*}
\end{small}

One can readily show that $\frac{15n^2+41n+60}{24} \geq 3n+3$ if $n \geq 3$.
Thus we completed the proof.
\end{proof}
\begin{theorem}
\label{th:n1a4}
$T(n,1;a,4)$ is true for any $n\ge1$ and $a \geq 3$.
\end{theorem}
\begin{proof}
The proof is by induction on $n$ and $a$.
Note that, since $a \geq 3$,
$T(1,1;a,4)$ is true by Theorem~\ref{th:CGG}.
We have also proved that
$T(n,1;3,4)$ and $T(n,1;4;4)$ are true for any $n\ge1$ (see Lemmas \ref{th:n134s} and \ref{th:n144s}).

Assume now that $n\ge2$ and $s\in\{\underline{s},\overline{s}\}$.
Let $s'$ and $\varepsilon$ be the quotient and remainder in the
division of $s(n+1)-5{n+a \choose a}$ by $n+1$ respectively.
Note that $(n,1;a-2,4;s-s'-\varepsilon)$
is superabundant by Lemma~\ref{th:supab} and
$s'\ge\varepsilon$ by Lemma \ref{s'big} (i).
Thus $T(n,1;a,4;s)$ can be reduced to $T(n-1,1;a,4;s')$, $T(n,1;a-1,4;s-s')$ and
$T(n,1;a-2,4;s-s'-\varepsilon)$. By induction hypotheses, these
statements are all true.
The statement $T(n,1;a,4;s)$ is therefore  true by Theorem~\ref{th:horace}.
\end{proof}
\begin{theorem}
\label{th:n1ab}
If $a,b \geq 3$, then $T(n,1;a,b)$ is true for any $n\ge1$.
\end{theorem}
\begin{proof}
The statement $T(1,1;a,b)$ is true by Theorem~\ref{th:CGG}, because
$a,b\ge3$. Suppose now that $n\ge2$. The proof is by induction on $b$.
Note that $T(n,1;a,3)$ is true by Theorem \ref{th:n1a2d}
and $T(n,1;a,4)$ is true by Theorem~\ref{th:n1a4} for any $a\ge3$. Thus
we may assume that $b\ge5$.
It is enough to prove $T(n,1;a,b;\underline{s})$ and $T(n,1;a,b;\overline{s})$.
Assume that $s\in\{\underline{s},\overline{s}\}$.
Let $s'$ and $\varepsilon$ be the quotient and remainder in the
division of $s(n+2)-{n+a \choose a}b$ by $n+1$.
By induction hypothesis, $T(n,1;a,b-1;s-s')$ and
$T(n,1;a,b-2;s-s'-\varepsilon)$ are true.  Additionally,
Lemma~\ref{th:supab} implies that $(n,1;a,b-2;s-s'-\varepsilon)$ is
superabundant, because $n\ge2$. Lemma~\ref{s'big}~(i) implies that
$s'\ge\varepsilon$, since $a,b\ge3$.
Note that $T(n,0;a,b;s')$ is true if $n \geq 2$ with only
four exceptions by the Alexander-Hirschowitz theorem \cite{AH}.
Thus the statement follows immediately from Theorem~\ref{th:horace} if  $(n,a,s') \not\in \{(2,4,5),(3,4,9),(4,3,7),(4,4,14)\}$.

Since $s'$ and $\varepsilon$ are the quotient and  remainder in the division by $n+1$ respectively, we have the following equality:
\begin{equation}\label{relation}
(n+1)s'+\varepsilon=s(n+2)-{n+a \choose a}b \ \   \mbox{for $0\le\varepsilon\le n$}.
\end{equation}
If $n=3$ and $a=4$, then $s=\underline{s}=\overline{s}=7(b+1)$, and thus $s'=8$ and $\varepsilon=3$.
So the above argument implies that $T(3,1;4,b;s)$, and hence $T(3,1;4,b)$, is true for every $b$. The same idea,  however, cannot be applied to $(n,a)=(2,5)$, $(4,3)$ and $(4,4)$.
Therefore, for each $(n,a) \in \{(2,4), (4,3), (4,4)\}$ and for each  $s$ such that
\[
\left(n,a,\frac{s(n+2)-{n+a \choose a}b-\varepsilon}{n+1} \right)
\]
falls into one of the above cases,
we  need to prove that $T(n,1;a,b; s)$ holds in a different way.

Let $t$ and $\delta$ be the quotient and remainder in the division of
$s(n+2)-{n+a-1 \choose a-1}(b+1)$ by $n+1$ respectively.
Note that $(n,1;a-2,b;s-t-\delta)$ is
superabundant by Lemma~\ref{th:supab}
and $t\ge\delta$ by Lemma \ref{s'big} (i).
So in order to apply Theorem
\ref{th:horace},  we need only to check that $T(n-1,1;a,b;t)$, $T(n,1;a-1,b;s-t)$ and
$T(n,1;a-2,b;s-t-\delta)$ are true.
Below we will consider the above-mentioned three cases separately.

\vspace{2mm}
\noindent
(i)  Let $(n,a,s')=(2,4,5)$. From \eqref{relation}, we have
$15+\varepsilon=4s-15b$ and $0\le\varepsilon\le 2$.
This implies that we can assume that $s$ is an integer of the form
$s=\frac{{15}(b+1)+\varepsilon}{4}$ for some $\varepsilon \in \{0,1,2\}$. It  suffices to prove the truth of $T(2,1;4,b;s)$ for such an $s$.

Let $t$ and $\delta$ be the quotient and remainder in the division of $4s-10(b+1)$ by $3$ respectively.
Note that $T(1,1;4,b;t)$ is true by Theorem \ref{th:CGG},
$T(2,1;3,b;s-t)$ is true for the first part of the proof.
Moreover $T(2,1;2,b;s-t-\delta)$ holds by Lemma \ref{abrescia} if $b$ is odd.
By Theorem \ref{th:n124s},  the statement is also true for $b=2k$ if
$s-t-\delta\ge 3(k+1)$.
Therefore, we need only to verify that  this inequality holds.
Assume that $b=2k$. Then $\varepsilon=1$ and $15b \equiv 0 \ (\mbox{mod  $4$})$,
which implies that $k$ is even. Set $k=2\ell$, so that $b=4\ell$ and $s=15\ell+4$.
Note that we may assume that $\ell\ge2$, because $b\ge5$.
Thus
$t=\lfloor\frac{4s-10(b+1)}{3}\rfloor=\lfloor\frac{20\ell+6}{3}\rfloor\le 7\ell+2$.
By definition,  $\delta\le 2$.
Thus we have $s-t-\delta \ge (15\ell+4)-(7\ell+2)-2 = 8\ell \ge 6\ell+3= 3(k+1)$.

\vspace{2mm}
\noindent
(ii)
Let $(n,a,s')=(4,3,7)$. From \eqref{relation}, we have
$35+\varepsilon=6s-35b$ and $0\le\varepsilon\le 4$.
This implies that we may assume that $s$ is an integer of the form
$s=\frac{{35}(b+1)+\varepsilon}{6}$ with $0\le \varepsilon\le 4$. Thus we need to prove $T(4,1;3,b;s)$ for such an $s$.

Let $t$ and $\delta$ be the quotient and remainder in the
division of $6s-15(b+1)$ by $5$. We have that
$T(3,1;3,b;t)$ holds by the first part of the proof. By Lemma \ref{abrescia}
$T(4,1;2,b;s-t)$ is true if $b$ is odd.
Additionally,  if $b=2k$ is even, $T(4,1;2,b;s-t)$
holds by Theorem \ref{th:n124s} if $s-t\le 5k$.
Thus it remains only to prove that this inequality holds.
Since we know that $\varepsilon\le4$, we need to take the following two cases into account:
\begin{enumerate}
\item[(a)]$\varepsilon=1$, $k=3h$, $b=6h$, $s=35h+6$;
\item[(b)]$\varepsilon=3$, $k=3h+1$, $b=6h+2$, $s=35h+18$.
\end{enumerate}
In case (a),
$t=\lfloor\frac{6s-15(b+1)}{5}\rfloor=
24h+\lfloor\frac{21}{5}\rfloor\ge 24h+4$. Thus
$s-t\le (35h+6)-(24h+4)=11h+2\le 15h=5k$.
In case (b),
$t=\lfloor\frac{6s-15(b+1)}{5}\rfloor \ge 24h+12$, and  hence
$s-t\le (35h+18)-(24h+12)=11h+6\le 5(3h+1)=5k$.

\vspace{2mm}
\noindent
(iii)
Let  $(n,a,s')=(4,4,14)$.
>From \eqref{relation}, we have $70+\varepsilon=6s-70b$ and
$0\le\varepsilon\le 4$.
This implies that $\varepsilon=2\varepsilon'$ is even, and
we may assume that $s$ is an integer of the form
$s=\frac{{35}(b+1)+\varepsilon'}{3}$ for $\varepsilon'=0,1,2$.

Let $t$ and $\delta$ be the quotient and remainder in the division of
$6s-35(b+1)$ by $5$.
We have already shown that  $T(3,1;4,b;t)$ and  $T(4,1;3,b;s-t)$ are true.
Thus we  only need to prove $T(4,1;2,b;s-t-\delta)$.
By Lemma \ref{abrescia} and Theorem \ref{th:n124s} , this statement holds  either if $b$ is odd or  if  $b=2k$ and $s-t-\delta\ge 5(k+1)$.
Thus all we have to do is to prove the above inequality holds.
To do so, we consider only the following possible three cases:
\begin{enumerate}
\item[(a)]
$k=3\ell$, $\varepsilon'=1$, $b=6\ell$, $s=70 \ell+12$,
\item[(b)]
$k=3\ell+1$, $\varepsilon'=0$, $b=6\ell+2$, $s=70 \ell+35$,
\item[(c)]
$k=3\ell+2$, $\varepsilon'=2$, $b=6\ell+4$, $s=70 \ell+59$.
\end{enumerate}
In all these three cases it is straightforward  to show that $s-t-\delta\ge 5(k+1)$ holds.
Thus we completed the proof.
\end{proof}
\begin{corollary}\label{coro3.14}
Let $n,a\ge1$, $b\ge3$, $\mathbf{n} = (n,1)$ and $\mathbf{a} = (a,b)$.
Then $X_{\mathbf{n},\mathbf{a}}$ is not defective except for $(n,a,b)=(n,2,2k)$.
\end{corollary}
\begin{proof}
In the previous theorem we proved the statement for $a,b\ge3$. So
we need only to consider the cases $a=1,2$.
Theorem \ref{n11d} implies that the statement is true if $a=1$.
By Remark \ref{rem-abrescia} and Theorem \ref{abrescia},
$X_{\mathbf{n},(2,b)}$ is defective if and only if $b$ is even.  Thus
we completed the proof.
\end{proof}
\begin{theorem}
\label{th:bound}
Suppose that $T(n,m;3,3)$, $T(n,m;3,4)$ and $T(n,m;4,4)$ are true
for any $n$ and $m$. Then
$T(m,n;a,b)$ is true for any $a, b\geq 3$.
\end{theorem}
\begin{proof}
We have already shown that $T(1,m;a,b)$ is true (see Theorem
\ref{th:n1ab}).  It follows from
Theorem~\ref{th:horace}, Lemma~\ref{th:supab}, and
Lemma~\ref{s'big}
that it is sufficient to prove that $T(n,m;3,b)$ and $T(n,m;4,b)$ are
true for every $b\ge3$.

We first prove that $T(n,m;3,b)$ is true for every $b\ge3$.
It has been already proved in Theorem~\ref{th:n1ab} that $T(n,1;3,b)$ is true.
We know by assumption that $T(n,m;3,3)$ and $T(n,m;3,4)$ are
true.
Thus the truth of $T(n,m;3,b)$ immediately follows from
Theorem~\ref{th:horace}.

We can analogously prove that $T(n,m;4,b)$ is true for every
$b\ge3$. Indeed, $T(n,1;4,b)$ holds by Theorem \ref{th:n1ab}, and
$T(n,m;4,3)$ and $T(n,m;4,4)$ are true by assumption.
Thus, by Theorem~\ref{th:horace}, $T(n,m;4,b)$ is also true.
\end{proof}
\section{Classification of $s$-defective Segre-Veronese
  varieties with $s \le 4$}
  \label{sec:classification}
This section is devoted to the classification of all the defective
$s^{\mathrm{th}}$ secant varieties of Segre-Veronese varieties with
$s \in \{2 , 3, 4\}$. Let $k \in \N$ and let $\mathbf{n}=(n_1,
\dots, n_k), \  \mathbf{a}=(a_1, \dots, a_k) \in (\Z_{\geq 0})^k
\setminus \{(0, \dots, 0)\}$.
The defective $s^{\mathrm{th}}$ secant varieties
of Segre varieties, i.e., Segre-Veronese varieties $X_{\mathbf{n},
\mathbf{a}}$ with $\mathbf{a} = (1^k)$,  has been completely
classified for such an $s$ in~\cite{AOP2} and for $k\ge3$. On the other
hand it is well known that $T(n_1,n_2;1,1;s)$ is false if and only if $2\le s\le \min\{n_1 n_2\}$.
We thus restrict our attention to the classification of defective
$s^{\mathrm{th}}$ secant varieties of Segre-Veronese varieties
$X_{\mathbf{n}, \mathbf{a}}$ for $\mathbf{a} > (1^k)$.
Let us first reformulate Lemma \ref{small-s} (ii) as follows:
\begin{lemma}\label{reformulation}
Suppose that $k \geq 2$ and that $s \in \{2,3,4\}$.
If the following are satisfied:
\begin{itemize}
\item[(i)] $\mathbf{m} \geq \mathbf{n}$;
\item[(ii)]  $\mathbf{a} > (1,1)$ if $k = 2$;
\item[(iii)]  $(\mathbf{n};\mathbf{a};s)$ is subabundant and
\item[(iv)]  $T(\mathbf{n};\mathbf{a};s)$ is true,
\end{itemize}
then $T(\mathbf{m};\mathbf{a};s)$ is also true.
\end{lemma}
\begin{proof}
If $\mathbf{n} = \mathbf{m}$, there is nothing to prove.
So we may assume that $\mathbf{n} < \mathbf{m}$. Let  $\Omega=\{ i \in
\{1, \dots, k\}\  | \ m_i > n_i  \  \mbox{and} \ a_i \geq 1\}$. Then
we showed that $\Omega \not= \emptyset$ in Lemma~\ref{small-s} (ii).
Lemma~\ref{small-s} (ii) also says that, in order to prove this lemma,
we only need to establish
\[
s \leq \min_{h \in \Omega}\left\{{n_h+a_h-1 \choose
a_h-1}\prod_{i\not=h} {n_i+a_i \choose a_i}\right\}.
\]
Without loss of generality, we may assume that
\[
{n_1+a_1-1 \choose a_1-1}\prod_{i=2}^k {n_i+a_i \choose a_i}=\min_{h
\in \Omega}\left\{{n_h+a_h-1 \choose a_h-1}\prod_{i\not=h} {n_i+a_i
\choose a_i}\right\}
\]
If $\mathbf{a} \ge (1^k)$ and $k \geq 3$, then
\begin{eqnarray*}
{n_1+a_1-1 \choose a_1-1} \prod_{i=2}^k {n_i+a_i \choose a_i}
&\geq &
 {n_1+a_1-1 \choose a_1-1} {n_2+a_2 \choose a_2} {n_3+a_3 \choose a_3} \\
& \geq &
{n_1+1-1 \choose 1-1} {n_2+1 \choose 1} {n_3+1 \choose 1}   \\
&\geq &  1 \, (n_2+1)(n_3+1) \\
&\geq  & 2\cdot 2 \\
& \geq & s.
\end{eqnarray*}
Suppose now that $k=2$. If $a_1 \geq 2$, then
\begin{eqnarray*}
{n_1+a_1-1 \choose a_1-1}  {n_2+a_2 \choose a_2}
&  \geq &  {n_1+2-1 \choose 2-1} {n_2+1 \choose 1} \\
& \geq &  (n_1+1)(n_2+1)  \\
& \geq & 2 \cdot 2 \\
& \geq & s.
\end{eqnarray*}
Similarly, if $a_1=1$ and if $a_{2} \geq 3$, then $ {n_1+a_1-1 \choose
a_1-1}  {n_2+a_2 \choose a_2}\ge 4$. Suppose now that $(a_1,
a_2)=(1,2)$. Analogously, we can immediately check that $ {n_1+a_1-1
\choose a_1-1}  {n_2+a_i \choose a_2}\ge 3$. Now assume that $s=4$
and note that $(n_1,1;1,2;4)$ is superabundant for every positive
integer $n_1$. Thus we may assume that $n_2 \geq 2$. Then it is
straightforward to see that also in this case $ {n_1+a_1-1 \choose a_1-1}
{n_2+a_2 \choose a_2}\ge 4=s$.
\end{proof}

For fixed $k\ge2$, $\mathbf{a} > (1,1,0,\ldots,0)$ and $s \in
\{2,3,4\}$, let us consider  the following partially ordered set:
\[
 M = \left\{\mathbf{n} \in \N^k \ | \
   \mbox{$(\mathbf{n};\mathbf{a};s)$ is subabundant}\right\}.
\]
Lemma \ref{reformulation} implies that, in order to prove that
$T(\mathbf{n};\mathbf{a};s)$ is true for every $\mathbf{n} \in M$,
it is enough to prove that $T(\mathbf{n};\mathbf{a};s)$ is true for
every minimal element of $M$ (there are only finitely many minimal
elements in $M$). In particular, if $(1^k; \mathbf{a}; s)$ is
subabundant and if $T(1^k; \mathbf{a}; s)$ is true, then
$T(\mathbf{n}; \mathbf{a}; s)$ is also true for every $\mathbf{n}
\in \N^k$.  In this case,   $T(\mathbf{m}; \mathbf{b}; s)$ is
also true for every $\mathbf{m} \in \N^{\ell}$ and for every
$\mathbf{b} \in \N^{\ell}$ with $\mathbf{b} \geq \mathbf{a}$ and
$\ell\ge k$, by Lemma \ref{small-s} (i) and Lemma \ref{reformulation}.

One can readily show that there are only finitely many superabundant
$k$-tuple $(\mathbf{n}; \mathbf{a}; s)$ if $s \in \{2,3,4\}$ except
for $(n,1;1,2;3)$, $(n,1;1,2;4)$ and $(n,1;1,3;4)$ with $n\geq 1$.
Since we have already proved, in Theorem \ref{n11d},
that $T(\mathbf{n};\mathbf{a};s)$ is true for each $(\mathbf{n};
\mathbf{a}; s) \in \{(n,1;1,2;3), (n,1;1,2;4), (n,1;1,3;4) \ | \  n
\geq 1\}$,  we only need to show the truth of a finite number of
statements to complete the classification of defective
$s^{\mathrm{th}}$ secant varieties of Segre-Veronese varieties for
the desired $s$.

In order to prove that $T(\mathbf{n}; \mathbf{a};
s)=S(\mathbf{n};\mathbf{a};s;0;0)$ is true for a given $(\mathbf{n};
\mathbf{a}; s)$, we apply Theorem~\ref{th:segrein}
that allows us to reduce it to proving
the truth of two statements of the forms
$S(\mathbf{n}';\mathbf{a};s';0;s'')$ and
$S(\mathbf{n}'';\mathbf{a};s'';0;s')$, where
$(\mathbf{n}';\mathbf{a};s';0;s'')$ and
$(\mathbf{n}'';\mathbf{a};s'';0;s')$ have the same abundancy. If the
truth of at least one of these statements, say
$S(\mathbf{n}';\mathbf{a};s';0;s'')$,  is not known yet,  then we
apply Theorem \ref{th:segrein} to $S(\mathbf{n}';\mathbf{a};s';0;s'')$. In order
to prove that  $T(\mathbf{n}; \mathbf{a}; s)$ is true, one must
repeat the same process over and over until  one achieves the
statements that are all known to be true. This procedure is
sometimes tedious to explicitly describe with words.  To avoid
tediousness,  we will represent this process by a tree diagram as
follows: Let $S(\mathbf{n};\mathbf{a};s;t;v)$ be a statement one
wishes to prove to be true.
Then the application of Theorem \ref{th:segrein} can be represented as the
following binary tree:
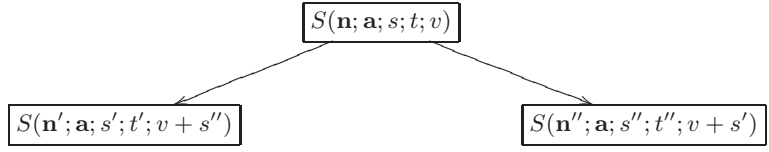
\begin{small}
\[
\xymatrix{
&*+[F]{S(\mathbf{n};\mathbf{a};s;t;v)}\ar[dl]  \ar[dr] \\
*+[F]{S(\mathbf{n}';\mathbf{a};s';t';v+s'')} & & *+[F]{S(\mathbf{n}'';\mathbf{a};s'';t'';v+s')}
}
\]
\end{small}
If the statements at the leaves of the tree are identical, we draw
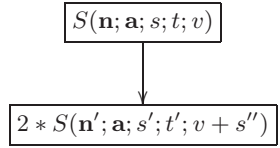
\begin{small}
\[
\xymatrix{
*+[F]{S(\mathbf{n};\mathbf{a};s;t;v)}\ar[d]  \\
*+[F]{2*S(\mathbf{n}';\mathbf{a};s';t';v+s'')}
}
\]
\end{small}

\noindent
instead of a usual binary tree.  In this case,
$(\mathbf{n};\mathbf{a};s;t;v)$ and
$(\mathbf{n}';\mathbf{a};s';t';v+s'')$ should have the same
abundancy.

The tree grows downward until one achieves only leaf nodes which are known to be true.
By Theorem \ref{th:segrein}, in order to prove that the statement at the root is
true, it suffices to show that the leaf statements are all true and
have the same abundancy.

\subsection{Case 1: $s=2$}
\begin{theorem}
\label{th:defective2} $T(\mathbf{n};\mathbf{a};2)$ is true with the
following exceptions:
\begin{itemize}
\item $k=1$, $n_1 \geq 2$ and $\mathbf{a}=2$;
\item $\mathbf{n}=(n_1,n_2)$ and $\mathbf{a}=(1,1)$ with $2 \leq n_1 \leq n_2$.
\end{itemize}
\end{theorem}
\begin{proof}
It is known by the Alexander-Hirschowitz theorem that if $k=1$, then
$T(\mathbf{n};\mathbf{a};2)$ fails if and only if $n_1 \geq 2$ and
$a_1=2$.  Thus we may assume that $k \geq 2$.

Suppose now that $\mathbf{a}=(1,2)$. Then
$(\mathbf{n};\mathbf{a};2)$ is subabundant for every $\mathbf{n} \in
\N^2$.  Since $T(1,1;1,2;2)$ is true by Theorem \ref{th:CGG},
it follows from Lemma~\ref{small-s} and Lemma~\ref{reformulation} that $T(\mathbf{n};\mathbf{a};2)$
are also true for all $\mathbf{n} \in \N^2$ and for all
$\mathbf{a}\ge(1,2)$.

Note that $(1^3;1^3;2)$ is equiabundant and $T(1^3;1^3;2)$ is true.
Hence $(\mathbf{n};\mathbf{a};2)$ are subabundant for all
$\mathbf{n}, \mathbf{a} \in \N^k$ with $k\geq 3$, and Lemma~\ref{small-s} and Lemma~\ref{reformulation}
imply that $T(\mathbf{n};\mathbf{a};2)$ is true, for any
$\mathbf{n}, \mathbf{a} > (1^3,0^{k-3})$.
\end{proof}
\subsection{Case 2: $s=3$}
\begin{proposition}
\label{th:(1^2,2;3)} The following statements are true:
\begin{itemize}
\item[(i)] $T(1^2,2;1^2,2;3)$;
\item[(ii)] $T(1,2,1;1^2,2;3)$.
\end{itemize}
\end{proposition}
\begin{proof}
Each statement can be reduced as in the following diagrams:

\vspace{1mm} \noindent (i)
{\small
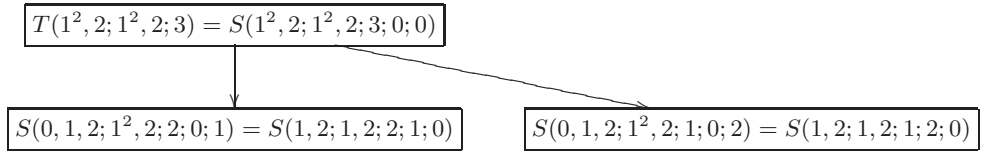
\[
\xymatrix{
*+[F]{T(1^2,2;1^2,2;3) =  S(1^2,2;1^2,2;3;0;0)} \ar[d] \ar[dr] \\
*+[F]{S(0,1,2;1^2,2;2;0;1)  = S(1,2;1,2;2;1;0)} &  *+[F]{S(0,1,2;1^2,2;1;0;2)  = S(1,2;1,2;1;2;0)}
}
\]
}

\vspace{1mm} \noindent (ii) {\small
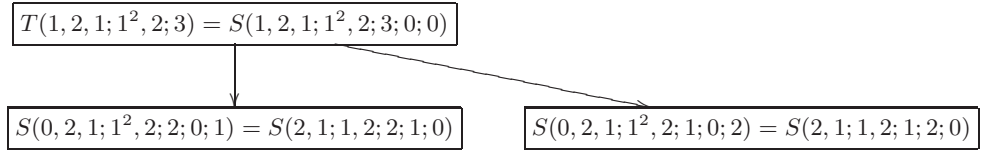
\[
\xymatrix{
*+[F]{T(1,2,1;1^2,2;3) =  S(1,2,1;1^2,2;3;0;0)} \ar[d] \ar[dr] \\
*+[F]{S(0,2,1;1^2,2;2;0;1)  = S(2,1;1,2;2;1;0)} &  *+[F]{S(0,2,1;1^2,2;1;0;2)  = S(2,1;1,2;1;2;0)}
}
\]
}

\vspace{1mm} \noindent
One can easily check that the following are all subabundant:
\[
\begin{array}{ccc}
(1^2,2;1^2,2;3),&
(1,2,1;1^2,2;3),  & (1,2;1,2;2;1;0), \\
(1,2;1,2;1;2;0), &
(2,1;1,2;2;1;0), & (2,1;1,2;1;2;0).
\end{array}
\]
By Theorem~\ref{th:defective2}, $S(1,2;1,2;2;0;0)$ and
$S(2,1;1,2;2;0;0)$ are true. Thus the following statements are also true
by Remark \ref{rem:remark} (v):
 \[
 S(1,2;1,2;2;1;0) \ \mbox{and} \ S(2,1;1,2;2;1;0).
 \]
 Likewise, $S(1,2;1,2;1;2;0)$ and
$S(2,1;1,2;1;2;0)$ are true, because $S(1,2;1,2;1;0;0)$ and
\linebreak
$S(2,1;1,2;1;0;0)$ are true. Thus we conclude that
\[
T(1^2,2;1^2,2;3) \ \mbox{and} \ T(1,2,1;1^2,2;3)
\]
are true.
\end{proof}
\begin{proposition}
\label{th:(1^4;1^3;2;3)} $T(1^4;1^3,2;3)$ is true.
\end{proposition}
\begin{proof}
We can reduce this statement as follows:
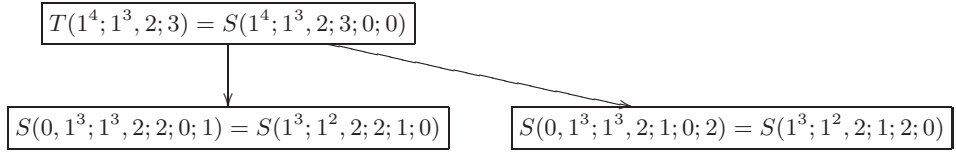
\[\small{
\xymatrix{
*+[F]{T(1^4;1^3,2;3) =  S(1^4;1^3,2;3;0;0)} \ar[d] \ar[dr] \\
*+[F]{S(0,1^3;1^3,2;2;0;1)  =  S(1^3;1^2,2;2;1;0)} & *+[F]{S(0,1^3;1^3,2;1;0;2)  =  S(1^3;1^2,2;1;2;0)}
}}
\]

\noindent Theorem~\ref{th:defective2}  implies that
$S(1^3;1^2,2;2;0;0)$ is true, from which it follows that
$S(1^3;1^2,2;2;1;0)$ is true.  Also, it is clear that
$S(1^3;1^2,2;1;0;0)$ is true, and so $S(1^3;1^2,2;1;2;0)$ is true
too. Since the following have the same abundancy:
\[
(1^3;1^2,2;2;1;0), \ (1^3;1^2,2;1;2;0) \ \mbox{and} \ (1^4;1^3,2;3),
\]
we can conclude that
$T(1^4;1^3,2;3)$ is true.
\end{proof}
\begin{theorem}
$T(\mathbf{n};\mathbf{a};3)$ is true with the following exceptions:
\begin{itemize}
\item $k=1$, $n_1 \geq 3$ and $\mathbf{a}=2$;
\item $\mathbf{n}=(n_1,n_2)$ with $3 \leq n_1 \leq n_2$ and $\mathbf{a}=(1,1)$;
\item $\mathbf{n}=(1,1)$ and $\mathbf{a}=(2,2)$;
\item $\mathbf{n}=(1,1,n)$ with $n \geq 3$ and $\mathbf{a}=(1,1,1)$;
\item $\mathbf{n}=(1,1,1)$ and $\mathbf{a}=(1,1,2)$;
\item $\mathbf{n}=(1^4)$ and $\mathbf{a}=(1^4)$.
\end{itemize}
\end{theorem}
\begin{proof}
Let $k=1$. Then from the theorem of Alexander and Hirschowitz it
follows that $T(\mathbf{n};\mathbf{a};3)$ is false if and only if
$n_1 \geq 3$ and $a_1=2$. Let us assume that $k \geq 2$.

The $11$-tuple $(1^5;1^5;3)$ is subabundant and $T(1^5;1^5;3)$ is
true (see~\cite{CGG4} for the proof). This means that if $k \geq 5$, then
$T(\mathbf{n};\mathbf{a};3)$ are true for all $\mathbf{n},
\mathbf{a} \in \N^k$.  Thus we may assume that $k \leq 4$.

Suppose that $k=4$. In~\cite{AOP2}, it was proved that if
$\mathbf{a}=(1^4)$, then there are no defective cases except for
$\mathbf{n}=(1^4)$.
We have proved in Proposition~\ref{th:(1^4;1^3;2;3)} that
$T(1^4;1^3,2;3)$ is true. This proves, by Lemma \ref{small-s} and
Lemma \ref{reformulation}, that
$T(\mathbf{n};\mathbf{a};3)$ is true for every $\mathbf{n}  \in
\N^4$ and $\mathbf{a}\geq (1^3,2)$. So the theorem holds if
$k=4$, and hence we may assume that $k \leq 3$.

Let $\mathbf{a}=(1^3)$. Then $T(\mathbf{n};\mathbf{a};3)$ is true
except for $\mathbf{n}=(1^2,n)$ with $n \geq 3$ (see~\cite{AOP2}). So
assume that $\mathbf{a} > (1^3)$. Since $T(1^3;1,2^2;3)$ is true
(see~\cite{BD}), it follows from Lemma \ref{small-s} and Lemma~\ref{reformulation} that if $\mathbf{a} \geq
(1,2^2)$, then  $T(\mathbf{n};\mathbf{a};3)$ is true for every
$\mathbf{n} \in \N^3$.
Note that $(1^3;1^2,3;3)$ is subabundant.
Additionally, the truth of $T(1^3;1^2,3;3)$ was proved by Bauer and
Draisma~\cite{BD}. Thus it remains only to show that
$T(\mathbf{n};1^2,2;3)$ is true except for $(\mathbf{n})=(1^3)$. It
is not hard to prove that $(\mathbf{n};1^2,2;3)$ is subabundant for
every $\mathbf{n} \in \N^3$. Since $T(1^3;1^2,2;3)$ is false, we
need to show that $T(1^2,2;1^2,2;3)$ and $T(1,2,1;1^2,2;3)$ are true
and we did it in Proposition \ref{th:(1^2,2;3)} . Thus we may now
assume that $k=2$.

Suppose that $\mathbf{a}=(1,1)$. It is known that
$T(\mathbf{n};\mathbf{a};3)$ holds if and only if  $3 \leq n_1 \leq
n_2$. Suppose that $\mathbf{a}=(1,2)$. Then
$(\mathbf{n};\mathbf{a};3)$ is subabundant if $\mathbf{n} \neq
(2,1),(1,1)$.  It was already proved in~\cite{BD} that
$T(1,2;1,2;3)$ is true. Thus $T(\mathbf{n};1,2;3)$ is true of
$\mathbf{n} \geq (1,2)$.   By Theorem \ref{n11d}, $T(n,1;1,2;3)$  with $n
\geq 1$ is true, and hence $T(\mathbf{n};1,2;3)$ is true for every
$\mathbf{n} \in \N^2$.

Let $\mathbf{a}=(1,3)$. Then $(\mathbf{n};\mathbf{a};3)$ is
superabundant if and only if $\mathbf{n}=(1,1)$. Note that
$T(1,1;1,3;3)$ is true, by Theorem \ref{n11d}. It was also proved that both
$T(1,2;1,3;3)$ and $T(2,1;1,3;3)$ are true (see~\cite{BD}),
which implies that $T(\mathbf{n};1,3;3)$ is true for any $\mathbf{n}
\in \N^2$.

Next consider $\mathbf{a}=(1,4)$. Clearly,
$(\mathbf{n};\mathbf{a};3)$ is subabundant for every $\mathbf{n} \in
\N^2$.  Since $T(1,1;1,4;3)$ is true by Theorem \ref{n11d},
$T(\mathbf{n};1,b;3)$ is also true for each $b \geq 4$ and
$\mathbf{n} \in \N^2$.

Let $\mathbf{a}=(2,2)$. Then $T(1,1;\mathbf{a};3)$ is known to be
false by Theorem \ref{th:CGG}.
Let $\mathbf{n}=(1,2)$. Then $(\mathbf{n};\mathbf{a};3)$ is
subabundant. Since $(\mathbf{n};1,1;3)$ is subabundant and since
$T(\mathbf{n};1,1;3)$ is true by Theorem \ref{n11d}, $T(\mathbf{n};\mathbf{a};3)$ is also
true. This also proves that $T(\mathbf{n};\mathbf{a};3)$ is true if
$\mathbf{n} > (1,1)$ and if $\mathbf{a} \geq (2,2)$.  Thus it
remains only to prove that $T(1^2;2,3;3)$ is true, because
$(1^2;2,3;3)$ is subabundant.
But $T(1^2;2,3;3)$ is true by Theorem \ref{th:CGG}.
Therefore we can conclude that
if $k=2$, then $T(\mathbf{n};\mathbf{a};3)$ fails if and only if
$\mathbf{n}=(n_1,n_2)$ with $3 \leq n_1 \leq n_2$ and
$\mathbf{a}=(1,1)$ or $\mathbf{n}=(1,1)$ and $\mathbf{a}=(2,2)$.
\end{proof}
\subsection{Case 3: $s=4$}
Let $\mathbf{n}=(n^2,1)$ with $n \geq 2$ and let $\mathbf{a}=(1^2,2)$.
Then $T(\mathbf{n};\mathbf{a};n+2)$ is known to be false by \cite[Corollary 5.5]{CC2}.
Here we give a different, but shorter proof of the same result:
\begin{proposition}
$T(n^2,1;1^2,2;n+2)$ is false for every $n \geq 2$.
\end{proposition}
\begin{proof}
Let $\mathbf{n}=(n^2,1)$ and let $\mathbf{a}=(1^2,2)$. The
defectivity of $\sigma_{n+2}(X_{\mathbf{n},\mathbf{a}})$ can be
proved by the existence of a certain rational normal curve in
$X_{2n+2}$ passing through generic $(n+2)$ points of
$X_{\mathbf{n},\mathbf{a}}$.

For each $i \in \{1,2,3\}$, let
$\pi_i$ the canonical projection from $\P^{\mathbf{n}}$ to the
$i^{\mathrm{th}}$ factor of $\P^{\mathbf{n}}$. Given generic points
$p_1, \dots, p_{n+2} \in \P^{\mathbf{n}}$, let $q_{i} = \pi_3(p_i)
\in \P^1$.
 Since any ordered subset of $n+2$ points in general position
 in $\P^n$ is projectively equivalent to the ordered set
 $\{\pi_i(p_1), \dots, \pi_i(p_{n+2})\}$
 for $i \in \{1,2\}$, there is a
 rational normal curve
$\nu_{n,i}: \P^1 \rightarrow C_n \subset \P^n$ of degree $n$
such that $\nu_{n,i} (q_j)=\pi_i(p_j)$ for all $j \in \{1, \dots,
n+2\}$. Let $\nu=(\nu_{n,1},\nu_{n,2}, \textrm{id})$ and let $C =
\nu(\P^1)$. Then $C$ passes through $p_1, \dots, p_{n+2}$. The image
of $C$ under the morphism given by $\cO(1^2,2)$ is  a rational
normal curve of degree $2n+2(=n+n+2\cdot1)$ in $\P^{2n+2}$.
Thus we have
\begin{small}
\begin{eqnarray*}
 \dim \sigma_{2n+2}(X_{\mathbf{n},\mathbf{a}})  &\leq  & 2n+2 +(n+2)(2n+1-1) \\
 & = & 2n^2+6n+2 \\
 & <  & (n+2)(2n+2)-1 \\
 & = & 2n^2+6n+3,
\end{eqnarray*}
\end{small}
and so  $ \sigma_{2n+2}(X_{\mathbf{n},\mathbf{a}})$ is defective.
\end{proof}
The above proposition only proves that if
$\mathbf{n}=(n^2,1)$ with $n \geq 2$ and $\mathbf{a}=(1^2,2)$, then
\linebreak
$\dim \sigma_{n+2}(X_{\mathbf{n},\mathbf{a}}) \leq 2n^2+6n+2$.
Below
we will show that  the equality actually holds.
\begin{proposition}
\label{th:(n,n,1)}
Let $\mathbf{n}=(n^2,1)$ with $n \geq 2$ and let
$\mathbf{a}=(1^2,2)$. Then
\[\dim
\sigma_{n+2}(X_{\mathbf{n},\mathbf{a}})=2n^2+6n+2.
\]
\end{proposition}
\begin{proof}
The statement $T(\mathbf{n};\mathbf{a};n+2)$ can be reduced to
$S(n-1,n,1;1^2,2;n+1;0;1)$ and \linebreak
$S(0,n,1;1^2,2;1;0;n+1)=S(n,1;1,2;1;n+1;0)$. Since
$(n,1;1,2;1;n+1;0)$ is subabundant and $S(n,1;1,2;1;0;0)$ is clearly
true, it follows that $S(n,1;1,2;1;n+1;0)$ is also true.
We can reduce $S(n-1,n,1;1^2,2;n+1;0;1)$  to
\[
S(n-2,n,1;1^2,2;n;0;2) \ \mbox{and} \ S(n,1;1,2;1;n+1;0).
\]
We continue in this manner until we reduce to
\[
S(0,n,1;1^2,2;2;0;n)=S(n,1;1,2;2;n;0) \ \mbox{and} \ n*S(n,1;1,2;1;n+1;0).
\]
Since $(n,1;1,2;2;n-1;0)$ is equiabundant and since
$S(n,1;1,2;2;0;0)$ is true by Proposition~\ref{th:defective2},
$S(n,1;1,2;2;n-1;0)$ is also true.
This proves the truth of \linebreak
$S(n,1;1,2;2;n;0)$. The truth of $S(n,1;1,2;2;n;0)$
implies that the linear subspace spanned by
$\sigma_2(X_{(n,1),(1,2)})$ and $n$ generic points coincides with
$\P^{3(n+1)-1}$; while the truth of $S(n,1;1,2;1;n+1;0)$ implies
that the linear subspace spanned by $\sigma_1(X_{(n,1),(1,2)})$ and
$n+1$ generic points has dimension $2n+2=(n+2)+(n+1)-1$. Therefore,
\begin{eqnarray*}
\dim  \sigma_{n+2}(X_{\mathbf{n},\mathbf{a}}) & \geq &  3n+3+n(2n+3)-1 \\
 & = & 2n^2+6n+2.
\end{eqnarray*}
Thus  we obtain $\dim
\sigma_{n+2}(X_{\mathbf{n},\mathbf{a}})=2n^2+6n+2$.
\end{proof}

\begin{remark}
More generally, let $\mathbf{n} = (n,n,1)$ with $n \geq 2$, let $\mathbf{a} = (1,1,2d)$ with $d \geq 1$.
As a particular case of \cite[Example 3.7]{AB2}, we know that $\sigma_s(X_{\mathbf{n},\mathbf{a}})$ is defective for any $d(n+1)+1 \leq s  \leq d(n+1)+n$.
If $\lfloor \frac{(2d+1)(n+1)}{2} \rfloor \leq s \leq nd+n+d$,
then $(\mathbf{n};\mathbf{a};s)$ is superabundant. Thus $\dim \sigma_s(X_{\mathbf{n},\mathbf{a}})  $ is expected to be $(n+1)^2(2d+1)-1$.
However,  there exist a form
$f_1$  of multi-degree $(1,0,d)$ and a form $f_2$ of multi-degree  $(0,1,d)$, both of which  vanish at given $s$ generic simple points.
Thus the form $f=f_1f_2$  of multi-degree $(1,1,2d)$ vanishes at the $s$ generic double points. So $\dim \sigma_s(X_{\mathbf{n},\mathbf{a}})  < (n+1)^2(2d+1)-1$, and hence
$\sigma_s(X_{\mathbf{n},\mathbf{a}})$ is defective.
It is worth mentioning that
these defective cases were first found  by Catalisano, Geramita and Gimigliano in~\cite{CGG2}.
When
$d(n+1)+1 \leq s  \leq \lfloor \frac{(2d+1)(n+1)}{2} \rfloor-1$, the proof of the defectivity of $\sigma_s(X_{\mathbf{n},\mathbf{a}})$ is not immediate and we refer to \cite{AB2} for more details.
On the other hand,  one can slightly modify the proof of Proposition~\ref{th:(n,n,1)} to show that  $\sigma_s(X_{\mathbf{n},\mathbf{a}})$ is not defective  if $s \leq d(n+1)$ or if $s \geq nd+n+d+1$. Thus we can conclude
that $\sigma_s(X_{\mathbf{n},\mathbf{a}})$ is defective exactly when  $d(n+1)+1\le s\le d(n+1)+n$.
\end{remark}

\begin{proposition}
\label{th:(1,n;1,d;4)} Let $k \in \N$ and let $\mathbf{n}=(n_2,
\dots, n_k) \in \N^{k-1}$ and let $\mathbf{a}=(a_2, \dots, a_k)
\in \N^{k-1}$ with either $\mathbf{a} > (1,1,0,\ldots,0)$ or
$\mathbf{a} > (2,0,\ldots,0)$. Then $T(1,\mathbf{n};1,\mathbf{a};4)$
is true.
\end{proposition}
\begin{proof}
Let $\mathbf{n}=(n_2, \dots, n_k)$ and let $\mathbf{a}=(a_2, \dots,
a_k)$. The statement can be reduced as in the following diagrams:
\begin{small}
\[
\xymatrix{
*+[F]{T(1,\mathbf{n};1,\mathbf{a};4) =  S(1,\mathbf{n};1,\mathbf{a};4;0;0)}\ar[d] \\
*+[F]{2*S(0,\mathbf{n};1,\mathbf{a};2;0;2)  =  2*S(\mathbf{n};\mathbf{a};2;2;0)}
}
\]
\end{small}

\noindent
Then $(1,\mathbf{n};1,\mathbf{a};4)$ and
$(\mathbf{n};\mathbf{a};2;2;0)$ must have the same abundancy. By
Theorem~\ref{th:defective2},
$S(\mathbf{n};\mathbf{a};2;0;0)=T(\mathbf{n};\mathbf{a};2)$ is true
because of the assumption on $\mathbf{a}$. Thus
$S(\mathbf{n};\mathbf{a};2;2;0)$ is true, from which the truth of
$T(1,n;1,d;4)$ follows.
\end{proof}
The following is an immediate consequence of
Proposition~\ref{th:(1,n;1,d;4)}:
\begin{corollary}
\label{th:6corollaries} The following statements are true:
\begin{itemize}
\item[(1)]  $T(1,2;1,b;4)$ are true for $b  \in \{3,4\}$;
 \item[(2)] $T(1,3,1;1^2,2;4)$;
 \item[(3)] $T(1^3;1^2,2;4)$;
\item[(4)] $T(1^3;1,2^2;4)$;
 \item[(5)] $T(1^2,2;1^2,2;4)$;
 \item[(6)] $T(1,2,1;1^2,2;4)$;
 \item[(7)] $T(1^3;1^2,3;4)$.
\end{itemize}
\end{corollary}
\begin{proposition}
\label{th:threefactor4} The following statements are all true:
\begin{itemize}
 \item[(i)] $T(2,3,1;1^2,2;4)$;
 \item[(ii)] $T(2^3;1^2,2;4)$;
\end{itemize}
\end{proposition}
\begin{proof}
Each statement in this theorem can be reduced as indicated in the
following diagram:

\vspace{1mm} \noindent (i)
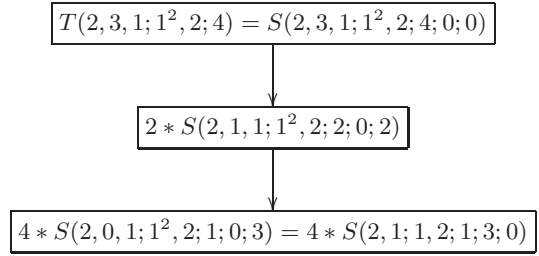
\begin{small}
\[
 \xymatrix{
*+[F]{T(2,3,1;1^2,2;4) =  S(2,3,1;1^2,2;4;0;0)}\ar[d] \\
*+[F]{2*S(2,1,1;1^2,2;2;0;2)}   \ar[d] \\
*+[F]{4*S(2,0,1;1^2,2;1;0;3) = 4*S(2,1;1,2;1;3;0)}
}
\]
\end{small}

\vspace{1mm} \noindent (ii)
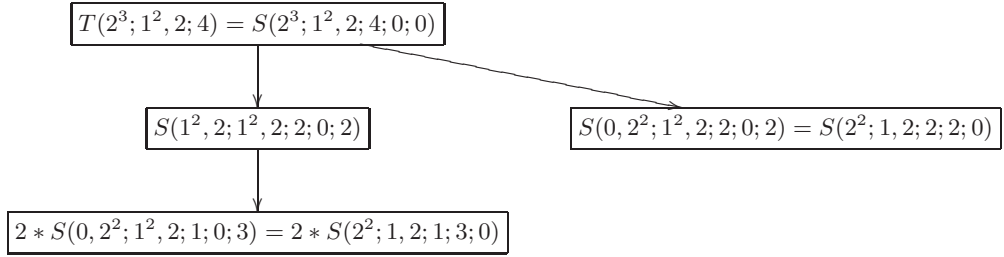
\begin{small}
\[
 \xymatrix{
  *+[F]{T(2^3;1^2,2;4) =  S(2^3;1^2,2;4;0;0)} \ar[d] \ar[dr] & \\
 *+[F]{S(1^2,2;1^2,2;2;0;2)}   \ar[d]  & *+[F]{S(0,2^2;1^2,2;2;0;2)=S(2^2;1,2;2;2;0)} \\
*+[F]{2*S(0,2^2;1^2,2;1;0;3) = 2*S(2^2;1,2;1;3;0)} &
}
\]
\end{small}

\noindent
Note that
the following are all subabundant:
\[
(2,3,1;1^2,2;4), \ (2,1;1,2;1;3;0), \ (2^3;1^2,2;4),
(2^2;1,2;2;2;0), \ (2^2;1,2;1;3;0).
\]
Since
$S(2^2;1,2;2;0;0)$  is true by Theorem~\ref{th:defective2}, so is
$S(2^2;1,2;2;2;0)$. Furthermore, both \linebreak
$S(2,1;1,2;1;0;0)$ and
$S(2^2;1,2;1;0;0)$ are true. Thus $S(2,1;1,2;1;3;0)$ and
$S(2^2;1,2;1;3;0)$ are also true. This means that all the statements
that appear at the leaf nodes in each tree are true. Thus also
$T(2,3,1;1^2,2;4)$ and $T(2^3;1^2,2;4)$ are true.
\end{proof}
\begin{proposition}\label{th:(1^4;1^3;2;4)}
$T(1^4;1^3,2;4)$ is true.
\end{proposition}
\begin{proof}
This statement can be reduced as follows:
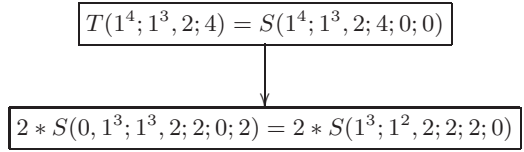
\begin{small}
\[
\xymatrix{
*+[F]{T(1^4;1^3,2;4) =  S(1^4;1^3,2;4;0;0)}\ar[d] \\
*+[F]{2*S(0,1^3;1^3,2;2;0;2)  =  2*S(1^3;1^2,2;2;2;0)}.
}
\]
\end{small}

\noindent
The classification of defective second secant varieties of
Segre-Veronese varieties (Theorem \ref{th:defective2}) implies that
$S(1^3;1^2,2;2;0;0)$ is true, from which the truth of
$S(1^3;1^2,2;2;2;0)$ follows. Since $(1^3;1^2,2;2;2;0)$ and
$(1^4;1^3,2;4)$ have the same abundancy, we
conclude that
$T(1^4;1^3,2;4)$ is true.
\end{proof}
\begin{theorem}
Let $k \in \N$ and let $\mathbf{n}, \mathbf{a} \in \N^k$. Then
$T(\mathbf{n};\mathbf{a};4)$ is false if and only if
$(\mathbf{n};\mathbf{a})$ falls into one of the following cases:
\begin{itemize}
\item $(\mathbf{n};\mathbf{a})=(n,2)$ with $n \geq 4$;
\item $\mathbf{n}=(n_1,n_2)$ with $4 \leq n_1 \leq n_2$ and $\mathbf{a}=(1,1)$;
\item $(\mathbf{n};\mathbf{a})=(1,2;2^2)$;
\item $(\mathbf{n};\mathbf{a})=(2^3;1^3)$;
\item $(\mathbf{n};\mathbf{a})=(1,2,n;1^3)$ with $n\geq4$;
\item $(\mathbf{n};\mathbf{a})=(2^2,1;1^2,2)$.
\end{itemize}
\end{theorem}
\begin{proof}
The Alexander-Hirschowitz theorem says that if $k=1$, then
$T(\mathbf{n};\mathbf{a};4)$ fails if and only if
$(\mathbf{n};\mathbf{a})=(n,2)$ with $n \geq 4$. Thus we may assume
that $k \geq 2$.

The $11$-tuple $(1^5;1^5;4)$ is subabundant and the truth of
$T(1^5;1^5;4)$ has been proved by Catalisano, Geramita and
Gimigliano (see~\cite{CGG4}). This means that if $k \geq 5$, then
$T(\mathbf{n};\mathbf{a};4)$ is true for all $\mathbf{n}, \mathbf{a}
\in \N^k$.  Thus we may assume that $k \leq 4$.

Suppose first that $k=4$. In~\cite{AOP2}, it was proved that there
are no defective cases if $\mathbf{a}=(1^4)$.  Note that
$(1^3,2;1^4;4)$ is equiabundant.  Thus it follows from Lemma \ref{small-s}
and Lemma \ref{reformulation} that $T(\mathbf{n};\mathbf{a};4)$ is true for every
$\mathbf{n} \in \N^k$ with $\mathbf{n} \geq (1^3,2)$ and for
every  $\mathbf{a} \in \N^k$. Moreover, we have already proved
that $T(1^4;1^3,2;4)$ is true (see
Proposition~\ref{th:(1^4;1^3;2;4)}). This proves that
$T(1^4;\mathbf{a};4)$ is true for every $\mathbf{a} \in \N^4$.
So the theorem holds if $k \geq 4$, and thus we may assume that $k
\leq 3$.

Suppose now that $k=3$.  In~\cite[Theorem 4.6]{AOP2},
$T(\mathbf{n};1^4;3)$ was proved to be true except for
$\mathbf{n}=(2^3)$ and $\mathbf{n}=(1,2,n)$ with $n \geq 4$.  So we
may assume that
$\mathbf{a} >(1^3)$.
In Corollary~\ref{th:6corollaries} it is proved that
$T(1^3;1,2^2;4)$ and $T(1^3;1^2,3;4)$ are true. Since
$(1^3;1,2^2;4)$ and $(1^3;1^2,3;4)$ are both subabundant, it follows
that $T(\mathbf{n};\mathbf{a};4)$ is true for every $\mathbf{n} \in
\N^3$ if $\mathbf{a} \geq (1^2,3)$ or $\mathbf{a} \geq (1,2^2)$.
This means that it remains only to prove the truth of
$T(\mathbf{n};1^2,2;4)$ for every $\mathbf{n} \in \N^3$ except
for $\mathbf{n}=(2^2,1)$. A $7$-tuple $(\mathbf{n};1^2,2;4)$ is not
equiabundant, but superabundant precisely when $\mathbf{n} = (1^3)$
and $(1,2,1)$. Additionally, we have proved that $T(2^2,1;1^2,2;4)$
is false in Proposition~\ref{th:(n,n,1)}. Thus all we need to do is  show that
$T(\mathbf{n};1^2,2;4)$ are true for all $\mathbf{n} \in \{(1^3),
(1,2,1), (1,3,1), (1^2,2), (2^3), (2,3,1)\}$. Those statements were,
however, proved to be true in Corollary~\ref{th:6corollaries} and
Proposition~\ref{th:threefactor4}.

Finally, assume that $k=2$.
In \cite{BD} it was already proved that $T(1^2;\mathbf{a};4)$ is
true for every $\mathbf{a}\in \N^2$ and that
$T(1,2;\mathbf{a};4)$ is true for every $\mathbf{a}\in \N^2$
except for $\mathbf{a}=(2,2)$.
Since the $5$-tuple $(1^2;2,3;4)$ is equiabundant and the statement
$T(1^2;2,3;4)$ is true, we conclude, by Lemma~\ref{small-s} and Lemma~\ref{reformulation},
that $T(\mathbf{n};\mathbf{a};4)$ is true for every $\mathbf{a} \geq
(2,3)$ and every $\mathbf{n}\in \N^2$. This means that to
complete the proof it is enough to prove the truth of $T(\mathbf{n};
\mathbf{a}; 4)$ for $\mathbf{a}=(1,d)$ with $d \geq 2$ and for
$\mathbf{a}=(2,2)$.

Assume first that $\mathbf{a}=(2,2)$. The $5$-tuple $(2^2;2^2;4)$ is
subabundant and the statement $T(2^2;2^2;4)$ is true by
Example \ref{segre2} and Remark \ref{rem:remark} (ii), because
the $5$-tuple $(2^2;2^2;5)$ is also subabundant.
So, by Lemma \ref{reformulation}, we conclude that $T(\mathbf{n};2^2;4)$ is true for every $\mathbf{n}\ge(2,2)$.

We now consider the case $\mathbf{a}=(1,d)$. Note that
$T(1^2;1,5;4)$ is true and $(1^2;1,5;4)$ is
  equiabundant.
Thus, by Lemma~\ref{small-s} and Lemma~\ref{reformulation},
$T(\mathbf{n};1,d;4)$ is also true for all $d
\geq 5$ and for any $\mathbf{n}\in\N^2$. If $d=4$, $(1^2;1,4;4)$
and $(2,1;1,4;4)$ are the only non-subabundant $5$-tuples. Thus
$T(\mathbf{n};1,4;4)$ is true for every $\mathbf{n}$, because the
truth of $T(\mathbf{n};1,4;4)$ was already proved to be true for
every $\mathbf{n} \in \{(1^2), (2,1), (3,1), (1,2)\}$.  Let $d=3$.
It is straightforward to prove that $(\mathbf{n};1,3;4)$ is not
subabundant if and only if
$\mathbf{n}=(m,1)$ with $m\geq 1$,
and $T(m,1;1,3;4)$ is true for every $m \geq 1$, by Theorem \ref{n11d}.
Furthermore, we proved in Corollary~\ref{th:(1,n;1,d;4)} that
$T(1,2;1,3;4)$ is true. This means that $T(\mathbf{n};1,3;4)$ holds
for every $\mathbf{n} \in \N^2$. Finally, suppose that $d=2$. It
is immediate to show that $(\mathbf{n};1,2;4)$ is not subabundant if
and only if $\mathbf{n}=(1,2)$, $\mathbf{n}=(2,2)$ or
$\mathbf{n}=(m,1)$ with $m \geq 1$.  This means that, in order to
prove the truth of $T(\mathbf{n};1,2;4)$ for every $\mathbf{n} \in
\N^2$, it is sufficient to show that $T(\mathbf{n};1,2;4)$ is
true for every $\mathbf{n} \in \{(1,2), (1,3), (2,2), (2,3), (m,1) \
\mbox{with $m \geq 1$} \}$. In~\cite{AB}, it was proved that
$T(1,n;1,2;4)$ are true for $n =\{2,3\}$ and that $T(2,3;1,2;4)$ is
true.  The truth of $T(2,2;1,2;4)$ was shown in~\cite{A}. Finally
$T(m,1;1,2;4)$ is true for every $m \geq 1$, by Theorem \ref{n11d}. Thus we
completed the proof.
\end{proof}
\section{Conjectures}
\label{sec:conjecture}
The main purpose of this section is to give a conjectural complete
list of defective  two-factor Segre-Veronese varieties. The first part
of this section is devoted to  collecting some results on defective
secant varieties of Segre-Veronese varieties.  To start with, we would
like to consider the so-called ``unbalanced" Segre-Veronese varieties.
\begin{definition}
Let $\mathbf{n}=(n_1, \dots, n_k) \in \N^k$ and let $\mathbf{a}=(a_1,
\dots, a_{k-1},1) \in \N^k$.
\begin{itemize}
\item
 $(\mathbf{n}; \mathbf{a})$ is said to be {\it  balanced} if $n_k \leq
 \prod_{i=1}^{k-1} {n_i+a_i \choose a_i}-\sum_{i=1}^{k-1} n_i$.
 \item
  $(\mathbf{n}; \mathbf{a})$ is said to be {\it  unbalanced} if $n_k
  \geq \prod_{i=1}^{k-1} {n_i+a_i \choose a_i}-\sum_{i=1}^{k-1}
  n_i+1$.
\end{itemize}
\end{definition}
The notion of ``unbalanced" was first introduced for Segre varieties
(see for example~\cite{CGG3} and~\cite{AOP2}). Then it was  extended to
Segre-Veronese varieties in~\cite{CGG1}.  The following theorem was
proved by Catalisano, Geramita and Gimigliano:
\begin{theorem}[\cite{CGG1}]
Let $\mathbf{n}=(n_1, \dots, n_k) \in \N^k$ and let $\mathbf{a}=(a_1,
\dots, a_{k-1},1) \in \N^k$.  Suppose that $(\mathbf{n}; \mathbf{a})$
is unbalanced. Then
$T(\mathbf{n};\mathbf{a};s)$ fails if and only if
\begin{eqnarray}
\label{eq:unbalanced}
 \prod_{i=1}^{k-1} {n_i+a_i \choose a_i}-\sum_{i=1}^{k-1} n_i<s<
 \min\left\{n_k+1, \  \prod_{i=1}^{k-1} {n_i+a_i \choose
     a_i}\right\}.
\end{eqnarray}
\end{theorem}
\begin{remark}
Let $\mathbf{n}$ and $\mathbf{a}$ be as given in the above theorem.
Then $X_{\mathbf{n},\mathbf{a}}$ is defective if and only if
Inequalities~(\ref{eq:unbalanced}) have an integer solution. Since
$(\mathbf{n}, \mathbf{a})$ is unbalanced,
if $n_k+1 \leq \prod_{i=1}^{k-1} {n_i+a_i \choose a_i}$, then
(\ref{eq:unbalanced}) must have at least one integer solution.

Suppose now that $n_k+1 > \prod_{i=1}^{k-1} {n_i+a_i \choose a_i}$,
then (\ref{eq:unbalanced}) have an integer solution if and only if
\[
 \prod_{i=1}^{k-1} {n_i+a_i \choose a_i} - \left[
  \prod_{i=1}^{k-1} {n_i+a_i \choose a_i}- \sum_{i=1}^{k-1} n_i
\right]
  = \sum_{i=1}^{k-1} n_i  > 1.
\]
This inequality holds unless $k=2$ and $n_1=1$.  Thus if
$(\mathbf{n},\mathbf{a})$ is unbalanced and if $(k,n_1)\not =(2,1)$,
then $X_{\mathbf{n},\mathbf{a}}$ is defective.
\end{remark}
Many other examples of defective secant varieties of  two-factor Segre-Veronese varieties have also been discovered by several authors. In Table  \ref{exceptionlist} below we provide the list of such defective secant varieties.
\begin{table}[ht]
\begin{center}
\begin{tabular}{|c|c|c|c|c|}
\hline
 & $\mathbf{n}$ & $\mathbf{a}$ & $s$ & References   \\
\hline
(1) & $(2,2k+1)$ & $(1,2)$ & $3k+2$  & \cite{O}\\
\hline
(2) & $(4,3)$ & $(1,2)$ & $6$  & \cite{CCh} \\
\hline
(3) & $(1,2)$ & $(1,3)$ & $5$ & \cite{DF}, \cite{CCh} \\
\hline
(4) & $(1,n)$ & $(2,2)$ & $n+2 \leq s \leq  2n+1$ & \cite{CGG2}, \cite{CGG1}, \cite{CC}\\
\hline
(5) & $(2,2)$ & $ (2,2)$ & $7$, $8$ & \cite{CGG2}, \cite{CGG1} \\
\hline
(6) & $(2,n)$ & $(2,2)$ & $\left\lfloor \frac{3n^2+9n+5}{n+3} \right\rfloor\leq s \leq 3n+2$ & \cite{CGG2}, \cite{B} \\
\hline
(7) & $(3,3)$ & $(2,2)$ & $14$, $15$ &  \cite{CGG2}, \cite{CGG1}  \\
\hline
(8) & $(3,4)$  & $(2,2)$ & $19$  & \cite{B}  \\
\hline
(9) & $(n,1)$ & $(2,2k)$ & $kn+k+1 \leq s \leq kn+k+n$ &  \cite{Abr}  \\
\hline
\end{tabular}
\end{center}
\caption{}\label{exceptionlist}
\end{table}
\begin{remark}
For an explanation of the cases where the degree is $(1,2)$ we refer
to \cite[Remark 5.1]{AB}. The defective cases of degree $(2,2)$ are
explained in \cite[Section 3]{CGG1}.
\end{remark}

We are now in position to state our conjecture:
\begin{conjecture}
\label{th:conjecture}
Let $\mathbf{n}=(m,n) \in \N^2$, let $\mathbf{a}=(a,b) \in \N^n$ and
let $X_{\mathbf{n}, \mathbf{a}}$ be the Segre-Veronese variety $\P^m
\times \P^n$ embedded by $\mathcal{O}_{ \P^m \times
  \P^n}(\mathbf{a})$. Then $X_{\mathbf{n}, \mathbf{a}}$ is defective
if and only if $(\mathbf{n}, \mathbf{a})$ falls into one of the
following cases:
\begin{itemize}
 \item[(a)] $(\mathbf{n}; \mathbf{a})=(m,n;a,1)$ is unbalanced and $m\ge2$.
 \item[(b)] $\mathbf{n}=(1,n)$ and $\mathbf{a}=(2k, 2)$ with $k \geq 1$.
 \item[(c)] $\mathbf{n}=(4,3)$, $(2,n)$ with $n$ odd and $\mathbf{a}=(1,2)$.
 \item[(d)] $\mathbf{n}=(1,2)$ and $\mathbf{a}=(1,3)$.
 \item[(e)] $\mathbf{n}=(2,2)$, $(3,3)$, $(3,4)$ and $\mathbf{a}=(2,2)$.
\end{itemize}
\end{conjecture}
Evidence for this conjecture was provided by the quoted results of
many authors. Further evidence in support of the conjecture was
obtained via computation.  Theorem~\ref{th:bound} suggests the
following little weaker conjecture:
 \begin{conjecture}
  Let $\mathbf{n}$, $\mathbf{a}$ and $X_{\mathbf{n}, \mathbf{a}}$ be as given in Conjecture~\ref{th:conjecture}. If $\mathbf{a} \geq (3,3)$, there are no defective two-factor Segre-Veronese varieties $X_{\mathbf{n}, \mathbf{a}}$ for all $\mathbf{n} \in \mathbb{N}^2$.
 \end{conjecture}

A substantial amount of effort has been made to complete the list of defective secant varieties of  two-factor Segre-Veronese varieties $X_{\mathbf{n},\mathbf{a}}$ for a given $(\mathbf{n}, \mathbf{a})$.  Below we list the cases that have been fully understood. Please refer to Table~\ref{exceptionlist} for the exceptions.

\begin{table}[ht]
\begin{center}
\begin{tabular}{|c|c|c|c|}
\hline
$\mathbf{n}$ & $\mathbf{a}$ & Exceptions & References   \\
\hline
$(1,n)$ & $(1,2)$ & None  & \cite{CCh}\\
\hline
$(2,n)$ & $(1,2)$ & (1)  & \cite{AB} \\
\hline
$(n,n-1)$ & $(1,2)$ & (2) & \cite{A} \\
\hline
$(n,n)$ & $(1,2)$ & None & \cite{A} \\
\hline
$(k,n)$ & $(1,k+1)$ &  None & \cite{CGG2} \\
\hline
$(1,2)$ & $(1,b)$ & (3)  & \cite{DF} \\
\hline
$(n,1)$ & $(1,b)$ &  None & \cite{ChiCil} \\
\hline
$(m,n)$ & $(1,b)$ with $b \geq 3$ and $(m+n+1)|{n+b \choose b}$ & None & \cite{BCC} \\
\hline
$(n,1)$ & $(2,b)$ & (9) & \cite{Abr} \\
\hline
$(n,1)$ & $(3,b)$ & (3) & \cite{Abr} \\
\hline
$(1,1)$ & $(a,b)$ & (9)  & \cite{CGG2}\\
\hline
$(n,1)$ & $(a,b)$ with $b\ge3$ & (9)  & Theorem \ref{th:n1abintro}  \\
\hline
\end{tabular}
\end{center}
\caption{}\label{classification}
\end{table}

\vspace{3mm}
\noindent
{\it Acknowledgements}.
We warmly thank Maria Virginia Catalisano for sharing with us her
knowledge in the field of secant varieties. We also thank  Giorgio Ottaviani
for many useful discussions and suggestions. Finally, we would like to thank the anonymous referee for helpful suggestions.



\begin{thebibliography}{aaa}
%
\bibitem{A} {H. Abo},
{\it On non-defectivity of certain Segre-Veronese varieties}, J. Symbolic Comput. {\bf 45} (2010), 1254--1269.
%
\bibitem{AB} {H. Abo \and M. C. Brambilla},
{\it Secant varieties of Segre-Veronese varieties $\P^m \times \P^n$ embedded
  by $\cO(1,2)$}, Experiment. Math. {\bf 18} (2009), no.~3, 369--384.
%
\bibitem{AB2}
 {H. Abo \and M. C. Brambilla},
{\it New examples of defective secant varieties of Segre-Veronese varieties}, arXiv:1101.3202, to appear in Collect. Math.
%
\bibitem{AOP2} {H. Abo, G. Ottaviani \and C. Peterson},
{\it Induction for secant varieties of Segre varieties},
Trans. Amer. Math. Soc. {\bf 361} (2009) 767-792.
%
\bibitem{AOP1} {H. Abo, G. Ottaviani \and C. Peterson},
{\it Non-defectivity of Grassmannians of planes},
to appear in J. Algebraic Geom.
%
\bibitem{Abr} {S.~Abrescia},
{\it About defectivity of certain Segre-Veronese varieties},
Canad. J. Math. {\bf 60}
(2008), no. 5, 961--974.
%
\bibitem{AH} {J. Alexander \and A. Hirschowitz},
{\it Polynomial interpolation in several variables}, J. Algebraic Geom.
{\bf 4} (1995), no.\ 2, 201--222.
%
\bibitem{Bal} {E.~Ballico},
{\it On the non-defectivity and non weak-defectivity of Segre-Veronese
  embeddings of products of projective spaces},
Port. Math. {\bf 63} (2006), no. 1, 101--111.
%
\bibitem{BD} {K.~Baur \and J.~Draisma},
 {\em Secant dimensions of low-dimensional homogeneous varieties},
Adv. Geom. {\bf 10} (2010), no. 1, 1--29.
%
\bibitem{BDG} {K.~Baur, J.~Draisma \and  W.~de Graaf},
{\em Secant dimensions of minimal orbits: computations and conjectures},
Experiment. Math. {\bf 16} (2007), no. 2, 239--250.
%
\bibitem{BCC}{A.~Bernardi, E.~Carlini \and M.~V.~Catalisano},
{\it Higher secant varieties of $\P^n \times \P^m$ embedded in bi-degree $(1,d)$}, arXiv:1004.2614.
%
%
\bibitem{B} {C.~Bocci},
{\it Special effect varieties in higher dimension},
Collect. Math. {\bf 56} (2005)
no. 3, 299--326.
%
\bibitem{BO} {M. C. Brambilla \and G. Ottaviani},
{\it On the Alexander-Hirschowitz theorem},
J. Pure and Applied Algebra {\bf 212} {2008}, no.\ 5, 1229--1251.
%
\bibitem{CC} {E. Carlini \and M. V. Catalisano},
{\it Existence results for rational normal curves},
J. Lond. Math. Soc. (2) {\bf 76} (2007), no.\ 1, 73--86.
%
\bibitem{CC2} {E. Carlini \and M. V. Catalisano},
{\it On rational normal curves in projective space},
J. Lond. Math. Soc. (2) {\bf 80} (2009), no.\ 1, 1--17.
%
\bibitem{CCh}
{E. Carlini \and J. Chipalkatti},
{\it On Waring's problem for several algebraic forms},
Comment. Math. Helv. {\bf 78}
(2003), no. 3, 494--517.
%
\bibitem{CGG3}
{M. V. Catalisano, A. V. Geramita \and A. Gimigliano},
{\it Ranks of tensors, secant varieties of Segre varieties and fat points},
Linear Algebra Appl. {\bf 355} (2002), 263--285.
%
\bibitem{CGG2}
{M. V. Catalisano, A. V. Geramita \and A. Gimigliano},
{\it Higher secant varieties of Segre-Veronese varieties},
Projective varieties with unexpected properties,
Walter de Gruyter GmbH \& Co. KG, Berlin, 2005, 81--107.
%
\bibitem{CGG1}
{M. V. Catalisano, A. V. Geramita \and A. Gimigliano},
{\it On the ideals of secant varieties to certain rational varieties},
J. Algebra {\bf 319} (2008), no. 5, 1913--1931.
%
\bibitem{CGG4}
{M. V. Catalisano, A. V. Geramita \and A. Gimigliano},
{\it Secant varieties of $\P^1 \times \cdots \times \P^1$ are not
  defective for $n \geq 5$},
to appear in J. Algebraic Geom.
%
\bibitem{Ch}
{K. A. Chandler},
{\it A brief proof of a maximal rank theorem for generic double points in projective space},
Trans. Amer. Math. Soc. {\bf 353} {2001}, no.\ 5, 1907--1920 (electronic).
%
\bibitem{ChiCil}
{L. Chiantini \and C. Ciliberto},
{\it The Grassmannians of secant varieties of curves are not defective},
Indag. Math., New Ser. {\bf 13} (2002), no.\ 1, 23-28.
%
\bibitem{DF}
{C. Dionisi \and C. Fontanari},
{\it Grassman defectivity \`a la Terracini},
Matematiche {\bf 56} (2001), no. 2, 245-255.
%
\bibitem{GS}
{D.~Grayson \and M.~Stillman},
{\it Macaulay 2, a software system for research in algebraic geometry},
Available at {\tt http://www.math.uiuc.edu/Macaulay2/}.
%
\bibitem{Har}
{J.~Harris}, {\it Algebraic geometry. A first course},
Graduate Texts in Mathematics {\bf 133} Springer-Verlag, Berlin, 1992.
%
\bibitem{L}
J.~M.~Landsberg,
{\it The geometry of tensors with applications},
in preparation.
\bibitem{O} {G.~Ottaviani},
{\it  Symplectic bundles on the plane, secant varieties and L\"uroth
  quartics revisited},
in Quaderni di Matematica, vol.21, editors G. Casnati, F. Catanese and
R. Notari, Vector Bundles and Low Codimensional Subvarieties: State of
the Art and Recent Developments, Aracne, 2008.
%
\end{thebibliography}
\end{document}